\documentclass[11pt]{article}

\usepackage{amssymb,latexsym}
\usepackage{pb-diagram}
\usepackage{amsmath,amscd}
\usepackage{citesort}
\usepackage{theorem} 
\usepackage{mathrsfs}

\title{\bf Gauge equivalence of Dirac structures and symplectic groupoids}

\author{Henrique Bursztyn\thanks{Current address: Department of Mathematics,
        University of Toronto, Toronto, Ontario, M5S 3G3. E-mail: henrique@math.toronto.edu}
        \\[0.1cm]
        Mathematical Sciences Research Institute\\
        1000 Centennial Drive, Berkeley\\
        CA, 94720 \\[0.2cm]and\\[0.2cm]
        Olga Radko \thanks{E-mail: radko@math.berkeley.edu} \\[0.1cm]
         Department of Mathematics\\
        University of California, Berkeley\\
        CA, 94720}

\date{September, 2002}

\newcommand{\id}         {{\mathrm {Id}}}

\newcommand{\graph}      {{\mathrm {graph}}}
\newcommand{\SP} [1]     {{\left\langle {{#1}} \right\rangle}}

\newcommand{\Dir}        {{\mathrm{Dir}}}
\newcommand{\Lag}        {{\mathrm{Lag}}}
\newcommand{\Bil}        {\mathrm{Bil}}
\newcommand{\OmegaL}     {\Omega_{\scriptscriptstyle{L}}}
\newcommand{\piL}        {\pi_{\scriptscriptstyle{L}}}
\newcommand{\piB}        {\pi_{\scriptscriptstyle{B}}}
\newcommand{\Lv}         {L_{\scriptscriptstyle{V}}}
\newcommand{\Lw}         {L_{\scriptscriptstyle{W}}}
\newcommand{\Lom}        {L_{\scriptscriptstyle{\Omega}}}
\newcommand{\Lpi}        {L_{\scriptstyle{\pi}}}
\newcommand{\tauB}       {\tau_{\scriptscriptstyle{B}}}
\newcommand{\tauphiB}    {\tau_{\scriptscriptstyle{\phi^*B}}}
\newcommand{\For}        {\mathcal{F}}
\newcommand{\Back}       {\mathcal{B}}
\newcommand{\pr}        {{\mathrm{pr}}}

\newtheorem{lemma} {Lemma} [section]
\newtheorem{proposition} [lemma] {Proposition}

\newtheorem{theorem} [lemma] {Theorem}
\newtheorem{corollary} [lemma] {Corollary}
\newtheorem{definition}[lemma] {Definition}

\theorembodyfont{\rm} 

\newtheorem{example}[lemma] {Example}

\newtheorem{remark}[lemma]{Remark}

\newenvironment{proof}{{\sc Proof:}}{{\hspace*{\fill} $\square$\\}}

\numberwithin{equation}{section}

\begin{document}

\maketitle

\begin{abstract}
We study gauge transformations
of Dirac structures and the relationship between gauge and Morita equivalences
of Poisson manifolds. We describe how the symplectic structure of a symplectic
groupoid is affected by a gauge transformation of the Poisson structure on its identity section, 
and prove that gauge-equivalent integrable Poisson structures are Morita equivalent.
As an example, we study certain  generic sets of   Poisson structures on Riemann
surfaces: we find complete gauge-equivalence invariants for such 
structures  which, on the $2$-sphere, yield a complete
invariant of Morita equivalence.
\end{abstract}

\section{Introduction}

\label{sec:intro}

Dirac structures were introduced in \cite{Cou90,CW88} to provide
a geometric framework for the study of constrained mechanical systems.
Examples of Dirac structures on a manifold $M$ include pre-symplectic
forms, Poisson structures and foliations; in general, a Dirac structure
determines a singular foliation on \( M \) whose leaves carry a 
pre-symplectic structure. 

The notion of gauge equivalence of Dirac
structures was introduced in \cite{SeWe01} motivated by the study
of the geometry of Poisson structures ``twisted'' by a closed \( 3 \)-form.
(See also \cite{KlSt01,Park}, where such structures were introduced in connection with Poisson-sigma
models.) 
There is a natural way to modify a Dirac structure on $M$ by a closed $2$-form
$B \in \Omega^2(M)$: one adds the pull-back of $B$ to the pre-symplectic form on each leaf
of the foliation, and the resulting pre-symplectic foliation determines
a new Dirac structure. This operation is called a gauge transformation, and it defines
an action of the additive group of closed $2$-forms on Dirac structures. Two Dirac structures
are called gauge equivalent if they lie in the same orbit of this action.
As discussed in \cite{SeWe01}, gauge-equivalent
Dirac structures share a lot of important properties; for instance,
their corresponding Lie algebroids are always isomorphic. 

Gauge transformations of Poisson structures also arise in  some
quantization problems. 
Let $(M,\pi)$ be a Poisson manifold, and let 
\( \mathrm{Def}(M,\pi ) \) denote its
moduli space of equivalence classes of star products \cite{BFFLS78}.
It follows from Kontsevich's formality theorem \cite{Kon97} that 
\( \mathrm{Def}(M,\pi ) \) is in bijection with equivalence
classes of formal Poisson structures on \( (M,\pi ) \).
As discussed in \cite{Bu2001},
the classification of deformation quantizations of  \( (M,\pi ) \) up to Morita
equivalence can be expressed in terms of the orbits of a canonical
action of the Picard group 
\( \mathrm{Pic}(M)\cong H^{2}(M,\mathbb Z) \) on \( \mathrm{Def}(M,\pi ) \).
Following Kontsevich's correspondence, 
there is a Poisson counterpart of algebraic Morita equivalence given by 
an $H^2(M,\mathbb{Z})$-action on formal Poisson structures, and the 
results in \cite{Bu2001,BuWa2001,JSW2001}
indicate that this action is given by
gauge transformations. In this context, a natural question is how gauge
equivalence relates to the geometric notion of Morita equivalence
of Poisson manifolds \cite{Xu91}. (See \cite{Land01} for other aspects
of the relationship between algebraic and geometric Morita equivalence.) 

In this paper we discuss the notion of gauge equivalence in the realm
of symplectic dual pairs and symplectic groupoids, and study the relationship
between the notions of gauge and Morita equivalence of Poisson structures. 

The paper is organized as follows. 

In Section \ref{sec:dirac} we recall some basic facts about Dirac
structures and their functorial properties, and establish an equivariance
property of Dirac maps with respect to gauge transformations. 

In Section \ref{sec:predual} we extend the usual notion of a symplectic
dual pair \cite{We83}  to deal with  Dirac manifolds. This
more general notion, called a pre-dual pair, arises naturally when
one considers a gauge transformation of a Poisson manifold which is
 part of a symplectic dual pair. We show that, under natural regularity
conditions,  ordinary symplectic
dual pairs are obtained as quotients of pre-dual pairs.

In Section \ref{sec:gaugegroup} we study the effect of gauge transformations
of an integrable Poisson manifold on the symplectic structure of its
symplectic groupoid. We show that if \( (G,\Omega ) \) is a symplectic
groupoid and \( G_{0} \) is its identity section, endowed with its
natural Poisson structure \( \pi  \), then a gauge transformation
of \( \pi  \) by \( B\in \Omega ^{2}(G_{0}) \) changes the symplectic
form on \( G \) by \[
\Omega \mapsto \Omega +\alpha ^{*}B-\beta ^{*}B,\]
 where \( \alpha  \) and \( \beta  \) are the source and target
maps, respectively. 

Using this result, we show in Section \ref{sec:morita} that two integrable
gauge-equivalent Poisson structures are Morita equivalent. The converse
is clearly not true, as Morita equivalent Poisson structures
need not have the same leaf decomposition. We show that, in fact,
the converse is not true even if we consider the possibility of gauge
equivalence up to a Poisson diffeomorphism. 

Finally, in Section \ref{sec:applic}, 
we study gauge and Morita equivalence
of a certain generic set of Poisson structures on a compact connected oriented
surface \( \Sigma  \). We consider the space \( \mathscr {G}_{n}(\Sigma ) \)
of Poisson structures vanishing linearly on \( n \) smooth disjoint
curves on \( \Sigma  \)  and show that two
Poisson structures \( \pi  \), \( \pi '\in \mathscr {G}_{n}(\Sigma ) \)
vanishing on the same curves are gauge equivalent if and only if they
have the same corresponding modular periods. (A classification of these structures up
to isomorphism was obtained in \cite{Radko01}.)
This result provides a  sufficient condition for Morita equivalence in
\( \mathscr {G}_{n}(\Sigma ) \); for the case \( \Sigma =S^{2} \),
we  construct a complete Morita-equivalence invariant for such  Poisson structures. 

\noindent{\bf Acknowledgments}: We would like to thank 
Robert Bryant, Marius Crainic, J.-P. Dufour, Rui L. Fernandes, Andr\'e Henriques, Nicholas Proudfoot and
Alan Weinstein for valuable discussions and comments. We also thank the referee for comments and
corrections.

\section{Dirac structures}

\label{sec:dirac}

\subsection{Linear Dirac structures}

\label{subsec:lineardirac}

Let \( V \) be a finite-dimensional real vector space. Let us equip
\( V\oplus V^{*} \) with the symmetric pairing \begin{equation}
\label{eq:sym}
\SP {(x,\omega ),(y,\mu )}=\frac{1}{2}(\omega (y)+\mu (x)).
\end{equation}
 A \textbf{Dirac structure} on \( V \) is a subspace \( L\subset V\oplus V^{*} \)
which is maximally isotropic with respect to \( \SP {\; ,\; } \) \footnote{%
This is equivalent to \( L \) being isotropic and dim \( L \) =
dim \( V \).
}. As we will see later, it is  useful to think of \( (V\oplus V^{*},\SP {\; ,\; }) \)
as an ``odd'' symplectic vector space (in the sense of super geometry), 
in such a way that its lagrangian subspaces
correspond to Dirac structures on $V$.
We denote the set of Dirac structures on \( V \) by \( \Dir (V) \). 

Throughout the text, we will identify bilinear forms \( \Omega :V\times V\longrightarrow \mathbb {R} \)
(resp. \( \pi :V^{*}\times V^{*}\longrightarrow \mathbb {R} \)) with
linear maps \( \tilde{\Omega }:V\longrightarrow V^{*} \) (resp. \( \tilde{\pi }:V^{*}\longrightarrow V \))
by \( \widetilde{\Omega }(u)(v)=\Omega (u,v) \) (resp. \( \tilde{\pi }(\nu )(\eta )=\pi (\eta ,\nu ) \)).
Whenever the context is clear, we will just write \( \Omega  \) for
\( \widetilde{\Omega } \) (resp. \( \pi  \) for \( \tilde{\pi } \)). 

\begin{example}
If \( \Omega  \)
is a skew-symmetric bilinear form on \( V \) (resp. \( \pi  \) is
a skew-symmetric bilinear form on \( V^{*} \), i.e., a bivector on
\( V \)), then \( L=\graph ({\Omega })\subset V\oplus V^{*} \) (resp.
\( L=\graph ({\pi }) \)) is a Dirac structure on \( V \). Conversely,
any Dirac structure \( L \) satisfying \( L\cap V^{*}=\{0\} \) (resp.
\( L\cap V=\{0\} \)) defines a skew-symmetric bilinear form (resp.
bivector) on \( V \). Hence pre-symplectic structures (i.e., skew-symmetric
bilinear forms) and Poisson structures (i.e., bivectors)
 on vector spaces are examples of linear Dirac structures.
\end{example}

Consider the natural projections \( \rho :V\oplus V^{*}\longrightarrow V \),
and \( \rho ^{*}:V\oplus V^{*}\longrightarrow V^{*} \), and let \( L \)
be a Dirac structure on \( V \). It is easy to check that \cite{Cou90} \begin{gather}
\rho(L)^{\circ }= L \cap V^{*}, \label{eq:charac1} \\
\rho^{*}(L)^{\circ }= L \cap V, \label{eq:charac2}
\end{gather} where \( ^{\circ } \) denotes the annihilator. 

\begin{proposition}\label{prop:dualint} A Dirac structure \( L \)
on \( V \) is equivalent to either of the following:
\begin{itemize}
\item[i)] A pair \( (R,\Omega ) \), where \( R\subseteq V \) is a subspace
and \( \Omega  \) is a skew-symmetric bilinear form on \( R \). 
\item[ii)] A pair \( (K,\pi ) \), where \( K\subseteq V \) is a subspace
and \( \pi  \) is a bivector on the quotient \( V/K \). 
\end{itemize}
Moreover, these correspondences are such that \( R = \rho (L) \)
and \( K=V\cap L=\ker \Omega  \). 
\end{proposition}

\begin{proof} Let \( L \) be a Dirac structure on \( V \). We set
\( R=\rho (L)\subseteq V \) and define the skew-symmetric bilinear
form \( \OmegaL :R\longrightarrow R^{*} \) by 
\begin{equation}
\label{eq:omegaL}
\OmegaL (x)\doteq \eta |_{{R}},
\end{equation}
 where \( \eta \in V^{*} \) is any element such that \( (x,\eta )\in L \).
The map \( \OmegaL  \) is well defined by (\ref{eq:charac1}), and
\( \ker \OmegaL =V\cap L \). Similarly, we can define a skew-symmetric
form \( \piL  \) on \( \rho ^{*}(L) \) and notice that 
\( (\rho ^{*}(L))^{*}=V/\rho ^{*}(L)^{\circ }=V/{V\cap L} \).
So, for \( K\doteq V\cap L \), we have a well-defined bivector \( \piL  \)
on \( V/K \) with \( \ker \piL =L\cap V^{*} \). 

Conversely, given a pair \( (R,\Omega ) \) as in \( i) \), we set
\begin{equation}
L\doteq \{(x,\eta )\; |\; x\in R,\, \eta \in V^{*}\,\mbox{ such that }\,
\eta |_{{R}}=\Omega (x)\}.
\end{equation}
 A simple dimension count shows that \( L \) is a Dirac structure
satisfying \( \rho (L)=R \) and \( \OmegaL =\Omega  \). A similar
construction holds for a pair \( (K,\pi ) \) as in \( ii \)). \end{proof}

For \( L\in \Dir (V) \), we denote the corresponding skew-symmetric
bilinear form on \( \rho (L) \) by \( \OmegaL  \) and the corresponding
bivector on \( V/V\cap L \) by \( \piL  \). 

\begin{example} Let \( L=\graph (\Omega ) \), where \( \Omega  \)
is a skew-symmetric bilinear form on \( V \). Then \( R=\rho (L)=V \),
\( \OmegaL =\Omega  \), \( K=\ker \Omega  \), and \( \piL  \) is
the Poisson bivector on the reduced space \( V/K \) corresponding
to the symplectic form induced by \( \Omega  \). 
\end{example}

\begin{example} Let \( L=\graph (\pi ) \), where \( \pi  \) is a
bivector on \( V \). Then \( R=\pi (V^{*}) \), and \( \OmegaL  \)
is the natural symplectic form induced by \( \pi  \). In this case,
\( K=\ker \OmegaL =\{0\} \), and \( \piL =\pi  \). 
\end{example}

\subsection{Functorial properties of linear Dirac structures}

\label{subsec:func}

Linear Dirac structures have nice functorial properties: they can
be both pushed forward and pulled back. This is a consequence of the
dual characterizations of Dirac structures (Prop. \ref{prop:dualint})
in terms of bilinear forms and bivectors. 

We can also describe the functorial properties of Dirac structures
in terms of Weinstein's symplectic category \cite{We82}, where the
objects are symplectic vector spaces  and the morphisms are canonical
relations. We recall here the main ideas.

Let \( E,F \) and \( H \) be symplectic vector spaces.
A \textbf{canonical relation} between \( E \) and \( F \) is a lagrangian
subspace \( L\subset E\times \overline{F} \), where \( \overline{F} \)
is the vector space \( F \) with symplectic form multiplied by \( -1 \).
If \( L_{1}\subset E\times \overline{F} \) and \( L_{2}\subset F\times \overline{H} \),
the usual composition of relations \begin{equation}
\label{eq:rel}
L_{1}\circ L_{2}=\{(x,y)\in E\times H\; |\; \exists z\in F\,
\mbox{ such that }\,(x,z)\in L_{1}\,\mbox{ and }\,(z,y)\in L_{2}\}
\end{equation}
defines a lagrangian subspace 
of \( E\times \overline{H} \). If \( L_{1}=\graph (f) \), \( L_{2}=\graph (g) \)
for symplectomorphisms \( f:F\longrightarrow E \) and \( g:H\longrightarrow F \),
then \( L_{1}\circ L_{2}=\graph (f\circ g) \). In general, the composition
of canonical relations defines a map 
\begin{equation}
\label{eq:comp}
\circ :\Lag (E\times \overline{F})\times \Lag (F\times \overline{H})\longrightarrow \Lag (E\times \overline{H}),
\end{equation}
 where \( \Lag (U) \) denotes the set of lagrangian subspaces of
a vector space \( U \). We remark that the same ideas work when symplectic vector
spaces are replaced by  vector spaces equipped with a nondegenerate symmetric bilinear form
with zero signature (thought of as ``odd'' symplectic spaces); in this case,
lagrangian (i.e., maximally isotropic) subspaces still have half the dimension of the total
space.

Let \( V \) and \( W \) be vector spaces, and let \( \phi :V\longrightarrow W \)
be a linear map. Let \( E=(V\oplus V^{*},\SP {\; ,\; }) \) and \( F=(W\oplus W^{*},\SP {\; ,\; }) \),
regarded as ``odd'' symplectic vector spaces. We define two canonical
relations associated to \( \phi  \): 
\begin{gather}
\For\phi = \{(\phi(x),\eta,x, \phi^{*}\eta) \;|\; x \in V, \, \eta \in W^{*}\} 
\in \Lag(F\times \overline{E}),
\label{eq:Fmap}\\
\Back\phi = \{(x,\phi^{*}\eta,\phi(x),\eta) \;|\; x \in V, \, \eta \in W^{*}\} 
\in \Lag(E\times \overline{F}).\label{eq:Bmap}
\end{gather} 
(As will become clear below, the letters $\For$ and $\Back$ stand for {\it forward} and
{\it backward}.)

Since we have the natural identifications \( \Dir (V)\cong \Lag (E\times \{0\}) \) and
\( \Dir (W)\cong \Lag (F\times \{0\}) \), the composition of relations
(\ref{eq:comp}) immediately induces maps 
\begin{gather}
\For\phi : \Dir(V) \longrightarrow \Dir(W) \label{eq:inducedF},\\
\Back\phi : \Dir(W) \longrightarrow \Dir(V) \label{eq:inducedB}.
\end{gather} 
Explicitly, for \( \Lv\in \textrm{Dir}(V),\, \Lw\in \textrm{Dir}(W) \)
we have 
\begin{gather}
\For\phi(\Lv)= \{(\phi(x),\eta)\;|\; x\in V, \eta \in W^{*}, \; (x,\phi^{*}\eta)\in \Lv\},\label{eq:FLv}\\
\Back\phi(\Lw)= 
\{(x,\phi^{*}\eta)\;|\; x\in V, \eta \in W^{*}, \; (\phi(x),\eta)\in \Lw\}.\label{eq:BLw}
\end{gather}

\begin{example} If \( \Lv =\graph (\pi ) \) for a bivector \( \pi  \)
on \( V \), then \( \For \phi (\Lv )=\graph (\phi _{*}\pi ) \).
Analo\-gous\-ly, if \( \Lw =\graph (\Omega ) \) for a skew-symmetric
bilinear form \( \Omega  \) in \( W \), then \( \Back \phi (\Lw )=\graph (\phi ^{*}\Omega ) \).
\end{example}

We observe that the maps \( \For \phi  \) and \( \Back \phi  \)
are not inverse to each other in general. A simple computation shows
that if \( \phi  \) is \( 1 \)-\( 1 \), then \( \Back \phi \circ \For \phi =\id  \),
and if \( \phi  \) is onto, then \( \For \phi \circ \Back \phi =\id  \). 

\begin{proposition}\label{prop:funcdual} Let \( \phi :V\longrightarrow W \)
be a linear map, and let \( \Lv \in \Dir (V) \) and \( \Lw \in \Dir (W) \). 

\begin{itemize}
\item[i)] If \( \For \phi (\Lv )=\Lw  \), then \( \ker \Omega _{{\Lw }}=\phi (\ker \Omega _{{\Lv }}) \)
and \( \pi _{{\Lw }}=\phi _{*}(\pi _{{\Lv }}) \). 
\item[ii)] If \( \Lv =\Back \phi (\Lw ) \), then \( \rho (\Lv )=\phi ^{-1}(\rho (\Lw )) \)
and \( \Omega _{{\Lv }}=\phi ^{*}(\Omega _{{\Lw }}) \). 
\end{itemize}
\end{proposition}

\begin{proof} Using (\ref{eq:FLv}), it is easy to check that if \( \Lw =\For \phi (\Lv ) \),
then \( \ker \Omega _{{\Lw }}=W\cap \Lw =\{\phi (x)\; |\; x\in V,\, (x,0)\in \Lv \} \).
Since \( V\cap \Lv =\ker \Omega _{{\Lv }} \), it follows that \( \ker \Omega _{{\Lw }}=\phi (\ker \Omega _{{\Lv }}) \). 

As \( \eta \in \phi (V\cap \Lv )^{\circ } \) implies that \( \phi ^{*}\eta \in V\cap \Lv  \),
we can define \( \phi _{*}\pi _{{\Lv }} \) on \( W/\phi (V\cap \Lv ) \)
by \( \phi _{*}\pi _{{\Lv }}(\eta )=\phi (\pi _{{\Lv }}(\phi ^{*}\eta )) \),
where \( \eta \in (W/\phi (V\cap \Lv ))^{*}\cong W^{*}/{\phi (V\cap \Lv )^{\circ }} \). 

By definition of \( \pi _{{\Lv }} \) (see Prop. \ref{prop:dualint}),
\( \phi (\pi _{{\Lv }}(\phi ^{*}\eta ))=\phi (x) \), where \( x \)
is such that \( (x,\phi ^{*}\eta )\in \Lv  \). On the other hand,
if \( y=\pi _{{\Lw }}(\eta ) \), then \( (y,\eta )\in \Lw =\For \phi (\Lv ) \),
which is the case if and only if \( y=\phi (x) \) and \( (x,\phi ^{*}\eta )\in L \).
Therefore \( \pi _{{\Lw }}=\phi _{*}(\pi _{{\Lv }}) \) and \( i \))
is proven. 

The proof of \( ii \)) is analogous and is left for the reader. 
\end{proof}

\begin{corollary}\label{cor:dircpoisson} Let \( (V,\pi _{1}) \)
and \( (W,\pi _{2}) \) be Poisson vector spaces, and let \( L_{\pi _{i}}=\graph (\pi _{i}) \),
\( i=1,2 \). A linear map \( \phi :V\longrightarrow W \) is Poisson
(i.e., \( \varphi _{*}\pi _{1}=\pi _{2} \)) if and only if  \( \For \phi (L_{\pi _{1}})=L_{\pi _{2}} \).
\end{corollary}

Similarly, if \( (V,\Omega _{1}) \)
and \( (W,\Omega _{2}) \) are pre-symplectic vector spaces, then a linear map
\( \phi :V\longrightarrow W \) satisfies $\phi^*\Omega_2 = \Omega_1$ if and only if
$\Back \phi (L_{\Omega_2}) = L_{\Omega_1}$, where $L_{\Omega_i}=\graph (\Omega_{i}),\, i=1,2$.

This motivates the following definition. 
\begin{definition}\label{def:dirac}
Let \( V \) and \( W \) be vector spaces endowed with Dirac structures
\( \Lv ,\Lw  \), respectively. A linear map \( \phi :V\longrightarrow W \)
is called \textbf{forward Dirac} if \( \For \phi (\Lv )=\Lw  \).
\end{definition}

A map satisfying the analogous definition for \( \Back \phi  \) is
called \textbf{backward Dirac}. We remark that these definitions are
not equivalent. In this paper we will only deal with forward Dirac
maps, and we will refer to them simply as \textbf{Dirac maps}. 

\begin{example}\label{ex:quotientmap} Let \( L\in \Dir (V) \). As
observed in Prop. \ref{prop:dualint}$(ii)$, 
\( V/{\ker \Omega _{L}} \)
has an induced Poisson structure \( \piL \). The projection 
\( \mathrm{pr}:V\longrightarrow V/{\ker \Omega _{L}} \) is a Dirac map. 
\end{example}

\begin{lemma}\label{lem:quotmap} 
Let \( \Lv \in \Dir (V) \) and
\( \Lw \in \Dir (W) \). If \( \Phi :V\longrightarrow W \) is a Dirac
map, then it naturally induces a Poisson map 
\( \phi :V/{\ker \Omega _{\Lv }}\longrightarrow W/{\ker \Omega _{\Lw }} \).
\end{lemma}

\begin{proof} Let \( \mathrm{pr}_{{V}}:V\longrightarrow V/\ker \Omega _{\Lv } \)
and \( \mathrm{pr}_{{W}}:W\longrightarrow W/\ker \Omega _{\Lw } \)
be the natural projections. By Prop. \ref{prop:dualint}(\( i \)),
the map \( \phi :V/\textrm{ker}\Omega _{L_{V}}\to W/\textrm{ker}\Omega _{L_{W}} \)
given by \( \phi (\mathrm{pr}_{{V}}(x))\doteq \mathrm{pr}_{{W}}(\Phi (x)) \)
is well defined. By the observation in Example \ref{ex:quotientmap},
we get \[
\For \phi (L_{\pi _{L_{V}}})=\For \phi (\For \mathrm{pr}_{{V}}(\Lv ))=\For \mathrm{pr}_{{W}}(\For \Phi (\Lv ))=\For \mathrm{pr}_{{W}}(\Lw )=L_{\pi _{L_{W}}},\]
 where \( L_{\pi _{L_{V}}}=\textrm{graph}(\pi _{L_{V}}),\, L_{\pi _{L_{W}}}=\textrm{graph}(\pi _{L_{W}}) \)
for the Poisson structures \( \pi _{\Lv} \) and \( \pi _{L_{W}} \)
on \( V/\textrm{ker}\Omega _{L_{V}} \) and \( W/\textrm{ker}\Omega _{L_{W}} \),
respectively. Therefore, \( \phi  \) is a Poisson map.
 \end{proof}

\subsection{Gauge equivalence of linear Dirac structures}

\label{subsec:lineargauge}

Let \( \Bil (W) \) be the additive group of skew-symmetric bilinear
forms on a vector space \( W \). Following Weinstein and Severa \cite{SeWe01},
we consider the action \( \tau :\Bil (W)\times \Dir (W)\longrightarrow \Dir (W) \)
by \textbf{gauge transformations}, 
\begin{equation}
\label{eq:tau}
\tauB (L)=\{(x,\eta +B(x))\; |\; (x,\eta )\in L\},\; \; B\in \Bil (W).
\end{equation}
 Alternatively, using Prop. \ref{prop:dualint}, two Dirac structures
\( L_{1},L_{2} \) on \( W \) are \textbf{gauge equivalent} if \begin{equation}
\label{eq:altertau}
\rho (L_{1})=\rho (L_{2})\; \mbox{ and }\; 
\Omega _{L_{1}}=\Omega _{L_{2}}+B|_{\rho (L_{2})},
\end{equation}
for $ B\in \Bil (W)$ and $L\in \Dir(W)$.

A linear map \( \phi :V\longrightarrow W \) of vector spaces induces
an action of \( \Bil (W) \) on \( \Dir (V) \), since \( \phi ^{*}\textrm{Bil}(W)\subseteq \textrm{Bil}(V) \). 

\begin{lemma}\label{lem:equiv} The map \( \For \phi :\Dir (V)\longrightarrow \Dir (W) \)
is \( \Bil (W) \)-equivariant. 
\end{lemma}

\begin{proof} We must show that \( \For \phi (\tauphiB (L))=\tauB (\For \phi (L)) \)
for all \( L\in \textrm{Dir}(V) \), \( B\in \textrm{Bil}(W) \). 

Since \( \tauphiB (L)=\{(x,\eta +(\phi ^{*}B)(x))\; |\; (x,\eta )\in L\} \),
we have \[
\For \phi (\tauphiB (L))=\{(\phi (x),\xi )\; |\; (x,\phi ^{*}\xi )\in \tauphiB (L)\}=\{(\phi (x),\xi )\; |\; (x,\phi ^{*}\xi -(\phi ^{*}B)(x))\in L\}.\]
 On the other hand, \[
\tauB (\For \phi (L))=\{(\phi (x),\zeta +B(\phi (x)))\; |\; (x,\phi ^{*}\zeta )\in L\}=\{(\phi (x),\xi )\; |\; (x,\phi ^{*}\xi -\phi ^{*}(B(\phi (x))))\in L\}.\]
 Since \( (\phi ^{*}B)(x)=\phi ^{*}(B(\phi (x))) \), the result follows.
\end{proof}

A similar equivariance property holds for the map \( \Back \phi  \). 

Let \( (V,\Omega ) \) be a symplectic vector space, and let \( (W,\pi ) \)
be a Poisson vector space. Let \( \Lom =\graph (\Omega ) \) and \( \Lpi =\graph (\pi ) \)
be the corresponding Dirac structures. 

\begin{lemma}\label{lem:sympoiss} 
Let \( \phi :(V,\Omega )\longrightarrow (W,\pi ) \)
be a linear Poisson map, and let \( B\in \textrm{Bil}(W) \). Then
the restricted map 
\( \phi :\ker (\Omega +\phi ^{*}B)\longrightarrow \ker \Omega _{\tauB (L_{\pi })} \)
is an isomorphism. 
In particular, the form ${\Omega + \phi^*B}$
is symplectic if and only if \( \tauB (\Lpi ) \) corresponds to a Poisson structure. 
\end{lemma}

\begin{proof} By Lemma \ref{lem:equiv}, \( \For \phi (\tauphiB (\Lom ))=\tauB (\For \phi (\Lom ))=\tauB (\Lpi ) \),
since \( \phi  \) is Poisson. Hence, by Prop. \ref{prop:funcdual}
\( (i) \), \begin{equation}
\label{eq:ker}
\ker \Omega _{\tauB (L_{\pi })}=\phi (\ker (\Omega +\phi ^{*}B)),
\end{equation}
 and \( \phi :\ker (\Omega +\phi ^{*}B)\longrightarrow \ker \Omega _{\tauB (L_{\pi })} \)
is onto. 

On the other hand, since \( \phi ^{*}B|_{{\ker {\phi }}}=0 \), it
follows that \[
\ker (\phi |_{\ker (\Omega +\phi ^{*}B)})=\ker \phi \cap \ker (\Omega +\phi ^{*}B)=\ker \phi \cap \ker \Omega =0.\]
 So \( \phi  \) is injective. \end{proof}

In order to study dual pairs, it will be useful to collect a few results
on pre-symplectic orthogonals in the linear case. 

\begin{lemma}\label{lem:orthog} Let \( L\in \Dir (V) \) and \( \phi :V\longrightarrow W \)
be linear. Then 
\[
(\ker \phi \cap \rho (L))^{\Omega _{L}}=
\{x\in V\; |\; \exists \eta \in W^{*}\,\mbox{ such that }\, (x,\phi ^{*}\eta )\in L\}.
\]
\end{lemma}

\begin{proof} If \( x\in (\ker \phi \cap \rho (L))^{\Omega _{L}} \),
then \( \Omega _{L}(x)(y)=0 \) for all \( y\in \ker \phi \cap \rho (L) \).
The form \( \eta \in (\phi (\rho (L)))^{*} \) given by \( \eta (\phi (y))=\Omega _{L}(x)(y) \)
is well defined, \( \phi ^{*}\eta =\Omega _{L}(x) \), and hence \( (x,\phi ^{*}\eta )\in L \).
On the other hand, suppose \( (x,\phi ^{*}\eta )\in L \) for \( \eta \in W^{*} \).
Then \( \Omega _{L}(x)=\phi ^{*}\eta |_{{\rho (L)}} \), and if \( y\in \ker \phi \cap \rho (L) \),
we have \( \Omega _{L}(x)(y)=\eta (\phi (y))=0 \). \end{proof}

Let \( L\in \Dir (V),L_{i}\in \Dir (W_{i}) \),
and let \( J_{i}:V\longrightarrow W_{i} \), \( i=1,2 \), be Dirac
maps. 

\begin{lemma}\label{lem:twomaps} Suppose that the following orthogonality
condition on the \( J_{i} \)-fibers holds: \begin{equation}
\label{eq:orthfibers}
(\ker {J_{1}}\cap \rho (L))^{\Omega _{L}}=\ker {J_{2}}\cap \rho (L).
\end{equation}
 Then for any \( B\in \Bil (W_{2}) \), \( J_{1}:(V,\tau _{{J_{2}^{*}B}}(L))\longrightarrow (W_{1},L_{1}) \)
is a Dirac map. \end{lemma}

\begin{proof} We must show that \( \For J_{1}(\tau _{{J_{2}^{*}B}}(L))=L_{1} \).
Since these subspaces have the same dimension, it suffices to prove
that \( \For J_{1}(L)=L_{1}\subseteq \For J_{1}(\tau _{{J_{2}^{*}B}}(L)) \).
Recall that \begin{eqnarray*}
\For J_{1}(\tau _{{J_{2}^{*}B}}(L)) & = & \{(J_{1}(x),\eta )\; |\; x\in V,\, \eta \in W_{1}^{*},\, (x,J_{1}^{*}\eta )\in \tau _{{J_{2}^{*}B}}(L)\}\\
 & = & \{(J_{1}(x),\eta )\; |\; x\in V,\, \eta \in W_{1}^{*},\, (x,J_{1}^{*}\eta -(J_{2}^{*}B)(x))\in L\},
\end{eqnarray*}
 and \[
\For J_{1}(L)=\{(J_{1}(x),\eta )\; |\; x\in V,\, \eta \in W_{1}^{*},\, (x,J_{1}^{*}\eta )\in L\}.\]
 Suppose \( (J_{1}(x),\eta )\in \For J_{1}(L) \). By Lemma \ref{lem:orthog},
\( x\in (\ker {J_{1}}\cap \rho (L))^{\Omega _{L}}=\ker {J_{2}}\cap \rho (L) \).
So \( (J_{2}^{*}B)(x)=0 \), and therefore \( (x,J_{1}^{*}\eta -(J_{2}^{*}B)(x))=(x,J_{1}^{*}\eta )\in L \).
This implies that \( (J_{1}(x),\eta )\in \For J_{1}(\tau _{{J_{2}^{*}B}}(L)) \),
and the result follows. \end{proof}

Consider a  symplectic vector space
\( (V,\Omega ) \), Poisson vector spaces \( (W_{i},\pi _{i}) \),
\( i=1,2 \), and linear Poisson maps \( J_{i}:V\longrightarrow W_{i} \).
We call the diagram
\[\begin{diagram}
\node{}\node{(V,\Omega)}\arrow{sw,t}{J_1}\arrow{se,t}{J_2}\node{}\\
\node{(W_1,\pi_1)}\node{}\node{(W_2,\pi_2)}
\end{diagram}\]
a {\bf linear dual pair} if
\( (\ker {J_{1}})^{\Omega }=\ker {J_{2}} \). Let us fix
such a linear dual pair,
and let \( \Lom =\graph (\Omega ) \) and \( L_{{\pi _{i}}}=\graph (\pi _{i}) \).
Before we prove the main result of this section, we need the following lemma.

\begin{lemma}\label{lem:orthog2} Let \( B\in \Bil (W_{2}) \), and
let \( {\Omega' }=\Omega +J_{2}^{*}B \). Then \( (\ker {J_{2}})^{{\Omega' }}=(\ker {J_{2}})^{\Omega }=\ker J_{1} \)
and \( (\ker {J_{1}})^{{\Omega' }}=(\ker {J_{1}})^{\Omega }+\ker {{\Omega' }} \).
\end{lemma}

\begin{proof} Since \( J_{2}^{*}B|_{{\ker J_{2}}}=0 \), it follows
that \( (\ker J_{2})^{{\Omega' }}=(\ker J_{2})^{\Omega }=\ker J_{1} \).
By taking \( {\Omega' } \)-orthogonals, we get \( (\ker J_{1})^{{\Omega' }}=\ker J_{2}+\ker {{\Omega' }} \).
\end{proof}

We can now prove the main result of this section.

\begin{theorem}\label{thm:ldualpair} 
Let \( B_{i}\in \Bil (W_{i}) \),
\( i=1,2 \), and let \( \widehat{\Omega }\doteq \Omega +J_{1}^{*}B_{1}+J_{2}^{*}B_{2} \).
Then
\begin{itemize}
\item[(i)] the maps \( J_{i}:(V,L_{\widehat{\Omega }})
\longrightarrow (W_{i},\tau _{{B_{i}}}(L_{{\pi _{i}}})) \), \( i=1,2 \), are Dirac;

\item[(ii)] the form \( \widehat{\Omega } \) is symplectic if and only if
\( \tau _{{B_{i}}}(L_{{\pi _{i}}}) \), \( i=1,2 \), are Poisson; 

\item[(iii)] \( (\ker {J_{1}})^{\widehat{\Omega }}=\ker {J_{2}}+\ker {\widehat{\Omega }} \). 
\end{itemize}
\end{theorem}

\begin{proof} By Lemma \ref{lem:equiv}, \( J_{2}:(V,L_{\Omega +J_{2}^{*}B_{2}})\longrightarrow (W_{2},\tau _{{B_{2}}}(L_{{\pi _{2}}})) \)
is Dirac, and, by Lemma \ref{lem:twomaps}, \( J_{1}:(V,L_{\Omega +J_{2}^{*}B_{2}})\longrightarrow (W_{1},\pi _{1}) \)
is also Dirac. Again by Lemma \ref{lem:equiv}, it follows that \( J_{1}:(V,L_{\widehat{\Omega }})\longrightarrow (W_{1},\tau _{{B_{1}}}(L_{{\pi _{1}}})) \)
is Dirac. By Lemma \ref{lem:orthog2}, \( (\ker J_{2})^{{\Omega +J_{2}^{*}B_{2}}}=(\ker J_{2})^{\Omega }=\ker J_{1} \).
So condition (\ref{eq:orthfibers}) is satisfied, and Lemma \ref{lem:twomaps}
implies that \( J_{2}:(V,L_{\widehat{\Omega }})\longrightarrow (W_{2},\tau _{{B_{2}}}(L_{{\pi _{2}}})) \)
is a Dirac map. This proves \( i) \). 

In order to prove \( ii) \), let us assume that \( \widehat{\Omega } \)
is symplectic. It follows from \( i) \) and Prop. \ref{prop:dualint}\( (i) \)
that \( \tau _{{B_{1}}}(L_{{\pi _{1}}}) \) and \( \tau _{{B_{2}}}(L_{{\pi _{2}}}) \)
are Poisson. Conversely, if \( \tau _{{B_{1}}}(L_{{\pi _{1}}}) \)
is Poisson, then \( \Omega +J_{1}^{*}B_{1} \) is symplectic, by Lemma
\ref{lem:sympoiss}. It follows, again from Lemma \ref{lem:sympoiss},
that if \( \tau _{{B_{2}}}(L_{{\pi _{2}}}) \) is also Poisson, then
\( \Omega +J_{1}^{*}B_{1}+J_{2}^{*}B_{2} \) is symplectic. 

We now prove \( iii) \). Clearly, \( (\ker J_{1})^{\widehat{\Omega }}=(\ker J_{1})^{{\Omega +J_{2}^{*}B_{2}}}=\ker J_{2}+\ker (\Omega +J_{2}^{*}B_{2}) \),
by Lemma \ref{lem:orthog2}. Thus \( \ker J_{2}+\ker {\widehat{\Omega }}\subseteq (\ker J_{1})^{\widehat{\Omega }} \).
On the other hand, again by Lemma \ref{lem:orthog2}, \( \ker (\Omega +J_{2}^{*}B_{2})\subseteq \ker J_{1} \).
So \( \ker {\widehat{\Omega }}\cap \ker J_{1}=\ker (\Omega +J_{2}^{*}B_{2}) \),
and therefore \( \ker (\Omega +J_{2}^{*}B_{2})\subseteq \ker {\widehat{\Omega }} \).
Hence \( \ker J_{2}+\ker (\Omega +J_{2}^{*}B_{2})\subseteq \ker J_{2}+\ker {\widehat{\Omega }} \),
and the result follows. 
\end{proof}

\subsection{Dirac structures on manifolds}

\label{subsec:diracman}

Let \( M \) be a smooth manifold. A \textbf{Dirac structure} on \( M \)
\cite{Cou90} is a subbundle \( L\subset TM\oplus T^{*}M \) which
determines linear Dirac structures pointwise and whose sections are
closed under the \textbf{Courant bracket} \( [\; \, ,\, \; ]:\Gamma (TM\oplus T^{*}M)\times \Gamma (TM\oplus T^{*}M)\longrightarrow \Gamma (TM\oplus T^{*}M) \),
\begin{equation}
\label{eq:courant}
((X,\omega ),(Y,\mu ))\mapsto ([X,Y],L_{X}\mu -L_{Y}\omega +\frac{1}{2}d(\omega (Y)-\mu (X
)).
\end{equation}

Let \( (M,L_{{M}}) \) and \( (N,L_{{N}}) \) be Dirac manifolds.
A smooth map \( \phi :M\longrightarrow N \) is a (forward) \textbf{ Dirac map}
if \( \For T_{x}\phi ((L_{{M}})_{x})=(L_{{N}})_{\phi (x)} \) for
all \( x\in M \). 

\begin{example} Let \( \Omega \in \Omega ^{2}(M) \) (resp. \( \pi \in \chi ^{2}(M) \)).
As discussed in Section \ref{subsec:lineardirac}, 
\( L=\graph (\Omega )\subset TM\oplus T^{*}M \)
(resp. \( L=\graph (\pi ) \)) defines a pointwise linear Dirac structure.
In this case, the extra condition involving the Courant bracket in
the definition of a Dirac manifold is equivalent to the integrability
condition \( d\Omega =0 \) (resp. \( [\pi ,\pi ]=0 \), where $[\,,\,]$ is the
Schouten bracket). Hence pre-symplectic
and Poisson structures on \( M \) are particular cases of Dirac structures.
\end{example}

The Courant bracket (\ref{eq:courant}) does not satisfy the Jacobi
identity in general. However, the Jacobi identity does hold when this
bracket is restricted to sections of a Dirac subbundle \( L\subset TM\oplus T^{*}M \),
and it defines a Lie algebroid structure on \( L \) with anchor
map \( \rho |_{{L}} \), where \( \rho :TM\oplus T^{*}M\longrightarrow TM \)
is the natural projection. When \( L=\graph (\pi ) \), where \( \pi  \)
is a Poisson structure, the Lie algebroid structure on \( L \) is
isomorphic to the natural Lie algebroid structure on \( T^{*}M \)
via the projection \( \rho ^{*}:TM\oplus T^{*}M\longrightarrow T^{*}M \). 

A Dirac manifold \( (M,L) \) carries a (singular) pre-symplectic
foliation: the leaves are the orbits of the corresponding Lie algebroid
and the leafwise pre-symplectic structure is defined as in Prop. \ref{prop:dualint}(\( i \)).
This foliation is symplectic if and only if \( L=\graph (\pi ) \)
for a Poisson structure \( \pi  \). 

As in the linear case, one can think of a Dirac structure \( L \)
on \( M \) as a ``Poisson structure on the leaf space \( M/\mathcal{K} \)'',
where \( \mathcal{K} \) is the characteristic foliation of the leafwise
pre-symplectic form \( \Omega _{L} \). More precisely, if \( (M,L) \)
is a Dirac manifold, we define its set of \textbf{admissible functions}
by 
\begin{equation}
\label{eq:admiss}
\mathcal{A}\doteq 
\{f\in C^{\infty }(M)\; |\; df|_{\ker (\Omega _{L})}=0\}=\{f\in C^{\infty }(M)\; |\; df\in \rho ^{*}(L)\}.
\end{equation}

If \( f\in \mathcal{A} \), there exists \( X\in \chi (M) \) such
that \( (X,df)\in L \); we call \( X \) a \textbf{hamiltonian vector
field} of \( f \) and denote it by \( X_{f} \). Note that hamiltonian
vector fields of admissible functions are defined up to vector fields
in \( \ker \Omega _{L} \). If \( f,g\in \mathcal{A} \), then the
bracket 
\begin{equation}
\{f,g\}\doteq \Omega _{L}(X_{f},X_{g})
\end{equation}
 is well defined and makes \( \mathcal{A} \) into a Poisson algebra.
Clearly, if \( L \) comes from a Poisson structure, then \( \mathcal{A} \)
is just the  Poisson algebra \( (C^{\infty }(M),\{\; ,\; \}) \).
If the characteristic foliation \( \mathcal{K} \) of \( \Omega _{L} \)
is \textbf{simple}, i.e., if \( M/\mathcal{K} \) is a smooth manifold
and \( \mathrm{pr}:M\longrightarrow M/\mathcal{K} \) is a submersion,
then \( \mathcal{A}\cong C^{\infty }(M/\mathcal{K}) \) and there
is a naturally defined Poisson structure \( \pi _{L} \) on \( M/\mathcal{K} \)
in such a way that \( \mathrm{pr} \) is a Dirac map (Example \ref{ex:quotientmap}). 

Gauge transformations of Dirac structures are defined analogously
to the linear case: if \( L \) is a Dirac structure on \( M \) and
\( B \) a closed \( 2 \)-form, we set 
\begin{equation}
\label{eq:gaugeman}
\tauB (L)=\{(X,\eta +B(X)),\; (X,\eta )\in L\}.
\end{equation}
 The closedness of \( B \) guarantees that \( \tauB (L) \) satisfies
the integrability condition with respect to the Courant bracket. As
in the linear case, two Dirac structures on \( M \) in the same \( \tau  \)-orbit
are called \textbf{gauge equivalent}. 

As observed in \cite[Sec.~3]{SeWe01}, for a Poisson structure \( \pi  \)
on \( M \), \( \tauB (\pi ) \) is Poisson if and only if the endomorphism
\( 1+B\pi :T^{*}M\longrightarrow T^{*}M \) is invertible. In this
case, 
$$ 
\tauB (\pi )=\pi (1+B\pi )^{-1}.$$
Gauge-equivalent Dirac structures share many properties: for instance,
they have the same leaf decomposition
(though the pre-symplectic forms on the leaves differ by the pullbacks of
\( B \) (see (\ref{eq:altertau}))) and their corresponding Lie algebroids
are isomorphic \cite{SeWe01}. In particular, gauge-equivalent Poisson
structures have isomorphic Poisson cohomology.

\section{Pre-dual pairs and reduction}

\label{sec:predual}
A \textbf{dual pair} \cite{We83}
consists of a symplectic manifold
\( (S,\Omega ) \), Poisson manifolds \( (M_{1},\pi _{1}) \), \( (M_{2},\pi _{2}) \),
and Poisson maps \( J_{i}:S\longrightarrow M_{i} \), \( i=1,2 \)
with symplectically orthogonal fibers, i.e.,
\begin{equation}
\label{eq:symporth}
\ker T_{x}J_{1}=(\ker T_{x}J_{2})^{\Omega }.
\end{equation}
 A dual pair is called \textbf{full} if the maps \( J_{1},J_{2} \)
are surjective submersions, and \textbf{complete} if these Poisson
maps are complete (see e.g. \cite{SilWein99}). 

If \( M_{1}\stackrel{J_{1}}{\leftarrow }S\stackrel{J_{2}}{\rightarrow }M_{2} \)
is a full dual pair, then 
\begin{gather}
\ker T_{x}J_{1} =\{X_{J_{2}^{*}f},\; f \in C^{\infty}(M_{2})\}, \mbox{ and}\label{eq:ham}\\
\{J_{1}^{*}(C^{\infty}(M_{1})),J_{2}^{*}(C^{\infty}(M_{2}))\}=0. \label{eq:comm}
\end{gather} 
In order to deal with Dirac structures, 
we generalize the notion of a dual pair as follows. 

\begin{definition}\label{def:predual} 
A \textbf{pre-dual pair}
is a pre-symplectic manifold \( (S,\Omega ) \), Dirac manifolds \( (M_{1},L_{1}) \),
\( (M_{2},L_{2}) \) and Dirac maps \( J_{i}:S\longrightarrow M_{i} \),
\( i=1,2 \), such that 
\begin{equation}
\label{eq:preorth}
(\ker T_{x}J_{1})^{\Omega }=\ker T_{x}J_{2}+\ker \Omega \; \mbox{ and }\; 
(\ker T_{x}J_{2})^{\Omega }=\ker T_{x}J_{1}+\ker \Omega .
\end{equation}
\end{definition}
As in the case of dual pairs, we represent pre-dual pairs by a diagram
\[\begin{diagram}
\node{}\node{(S,\Omega)}\arrow{sw,t}{J_1}\arrow{se,t}{J_2}\node{}\\
\node{(M_1,L_1)}\node{}\node{(M_2,L_2).}
\end{diagram}\]

As before, a pre-dual pair will be called \textbf{full} if each \( J_{i} \)
is a surjective submersion. 

\begin{proposition}\label{prop:preimpliesdual} 
Consider a pre-dual
pair 
\( M_{1}\stackrel{J_{1}}{\leftarrow }S\stackrel{J_{2}}{\rightarrow }M_{2} \).
If \( S \) is symplectic, then \( M_{1} \) and \( M_{2} \) are
automatically Poisson, and we have a dual pair in the usual sense.
\end{proposition}

\begin{proof} If \( \ker (\Omega )=0 \), then 
\( \ker (\Omega _{L_{i}})=TJ_{i}(\ker \Omega )=0 \),
and hence \( L_{i} \) is automatically Poisson. The maps \( J_{i} \),$i=1,2$,
are Dirac, and hence Poisson, and the orthogonality property (\ref{eq:preorth})
reduces to (\ref{eq:symporth}). 
\end{proof}

We now fix a full pre-dual pair 
\( M_{1}\stackrel{J_{1}}{\leftarrow }S\stackrel{J_{2}}{\rightarrow }M_{2} \). 

\begin{proposition}\label{prop:preham} We have the following generalization
of (\ref{eq:ham}): \[
\ker T_{x}J_{1}+\ker \Omega =
(\ker T_{x}J_{2})^{\Omega }=\{X_{J_{2}^{*}f}\; |\; f\in C^{\infty }(M_{2})\}.\]
\end{proposition}

\begin{proof} Clearly, if \( v\in \ker T_{x}J_{2} \), \( \Omega (X_{J_{2}^{*}f},v)=d_{x}(J_{2}^{*}f)(v)=df(T_{x}J_{2}(v))=0 \).
On the other hand, if \( u\in (\ker T_{x}J_{2})^{\Omega } \), \( \eta =\Omega _{x}(u)\in T^{*}_{x}S \)
vanishes along \( \ker T_{x}J_{2} \). Since \( J_{2} \) is a submersion,
in a neighborhood of \( x \) we can choose an exact \( 1 \)-form
\( dg \) with \( d_{x}g=\eta  \) and such that \( dg \) vanishes
on the distribution \( \ker TJ_{2} \). Hence \( g \) is constant
along the \( J_{2} \)-fibers and can be written in the form \( J_{2}^{*}f \)
for \( f\in C^{\infty }(M_{2}) \). Around \( x \), \( u=X_{J_{2}^{*}f} \).
\end{proof}

The next result generalizes property (\ref{eq:comm}). 
\begin{proposition}\label{prop:preorth}
Let \( \mathcal{A} \) be the algebra of admissible functions on \( S \),
and let \( \mathcal{A}_{i} \) be the algebra of admissible functions
on \( M_{i} \). Then 
\begin{itemize}
\item[i)] \( J_{i}^{*}\mathcal{A}_{i}\subseteq \mathcal{A} \), 
\item[ii)] \( \{J_{1}^{*}\mathcal{A}_{1},J_{2}^{*}\mathcal{A}_{2}\}=0 \), 
\item[iii)]\( \{J_{i}^{*}f,J_{i}^{*}g\}=J_{i}^{*}\{f,g\} \), for \( f,g\in \mathcal{A}_{i} \). 
\end{itemize}
\end{proposition}

\begin{proof} By Prop. \ref{prop:funcdual}, \( TJ_{i}(X)\in \ker \Omega _{L_{i}} \).
So, if \( X\in \ker \Omega  \) and \( f\in \mathcal{A}_{i} \), then
\( d(J_{i}^{*}f)(X)=df(TJ_{i}(X))=0 \), and \( i) \) is proven. 

For \( ii) \), note that if \( f_{i}\in \mathcal{A}_{i} \), then
\( \{J_{1}^{*}f_{1},J_{2}^{*}f_{2}\}=\Omega (X_{J_{1}^{*}f_{1}},X_{J_{2}^{*}f_{2}})=0 \)
by Prop. \ref{prop:preorth}. 

Finally, since \( J_{i} \) is a Dirac map, \( (X,J_{i}^{*}dg)\in L \)
implies that \( (T_{x}J_{i}(X),dg)\in \For T_{x}J_{i}(L)=L_{i} \).
This means that \( T_{x}J_{i}(X_{{J_{i}^{*}g}})=X_{g} \). Now, if
\( f,g\in \mathcal{A}_{i} \), then \[
\{J_{i}^{*}f,J_{i}^{*}g\}(x)=d_{x}(J_{i}^{*}f)(X_{{J_{i}^{*}g}})=
d_{J_i(x)}f(T_{x}J_{i}(X_{{J_{i}^{*}g}}))=d_{{J_{i}(x)}}f(X_{g})=\{f,g\}(J_{i}(x)),\]
 and the result follows. \end{proof}

Suppose that the characteristic foliations of \( \Omega  \) and \( \Omega _{L_{i}} \)
, \( \mathcal{K} \) and \( \mathcal{K}_{i} \), respectively, are
simple. Recall that \( S/\mathcal{K} \) has an induced symplectic
structure $\Omega_{\mathrm{red}}$ , while \( M_{i}/\mathcal{K}_{i} \), \( i=1,2 \), carry
induced Poisson structures $\pi_i$ so that the projections \( \mathrm{pr}:S\longrightarrow S/\mathcal{K} \)
and \( \mathrm{pr}_{i}:S\longrightarrow M_{i}/\mathcal{K}_{i} \)
are Dirac maps. 

\begin{theorem}\label{thm:reduction} 
Let \( M_{1}\stackrel{J_{1}}{\leftarrow }S\stackrel{J_{2}}{\rightarrow }M_{2} \)
be a full pre-dual pair. The maps \( J_{i}:S\longrightarrow M_{i} \)
induce surjective submersions on the quotient, 
\( j_{i}:S/\mathcal{K}\longrightarrow M_{i}/\mathcal{K}_{i} \),
in such a way that 
\[\begin{diagram}
\node{}\node{(S/\mathcal{K},\Omega_{\mathrm{red}})}\arrow{sw,t}{j_1}\arrow{se,t}{j_2}\node{}\\
\node{(M_1/\mathcal{K}_1,\pi_1)}\node{}\node{(M_2/\mathcal{K}_2,\pi_2)}
\end{diagram}\]
is a full dual pair. 
\end{theorem}

\begin{proof} By Prop. \ref{prop:funcdual}\( (i) \), \( TJ_{i}(\ker \Omega )=\ker \Omega _{L_{i}} \).
Hence if \( x \) and \( y \) belong to the same leaf of \( \mathcal{K} \),
then \( J_{i}(x) \) and \( J_{i}(y) \) belong to the same leaf of
\( \mathcal{K}_{i} \). So the map \( j_{i}:S/\mathcal{K}\longrightarrow M_{i}/\mathcal{K}_{i} \),
\( j_{i}(\mathrm{pr}(x))=\mathrm{pr}_{i}(J_{i}(x)) \) is well defined
and is a surjective submersion. It follows from Lemma \ref{lem:quotmap}
that \( j_{i} \) is a Poisson map. 
Finally,
a simple computation shows that \( (\ker Tj_{i})^{\Omega _{\mathrm{red}}}=\mathrm{pr}((\ker TJ_{i})^{\Omega }) \).
Hence (\ref{eq:preorth}) implies that \( (\ker Tj_{1})^{\Omega _{\mathrm{red}}}=\ker Tj_{2} \).
\end{proof}

We obtain examples of pre-dual pairs from gauge transformations of Poisson
structures in dual pairs as follows. 

\begin{theorem}\label{thm:gaugedual} Let \( (M_{1},\pi _{1})\stackrel{J_{1}}{\leftarrow }(S,\Omega )\stackrel{J_{2}}{\rightarrow }(M_{2},\pi _{2}) \)
be a full dual pair, and let \( B_{i} \) be a closed \( 2 \)-form
on \( M_{i} \), \( i=1,2 \). Let \( \widehat{\Omega }=\Omega +J_{1}^{*}B_{1}+J_{2}^{*}B_{2} \).
Then 
\[\begin{diagram}
\node{}\node{(S,\widehat{\Omega })}\arrow{sw,t}{J_1}\arrow{se,t}{J_2}\node{}\\
\node{(M_1,\tau _{{B_{1}}}(L_{{\pi _{1}}}))}\node{}\node{(M_2,\tau _{{B_{2}}}(L_{{\pi _{2}}}))}
\end{diagram}\]
is a full pre-dual pair. Moreover, \( \widehat{\Omega } \) is symplectic
if and only if \( \tau _{{B_{i}}}(L_{{\pi _{2}}}) \), \( i=1,2 \),
are Poisson, in which case they form a dual pair. \end{theorem}

\begin{proof} The result is a consequence of Theorem \ref{thm:ldualpair}.
\end{proof}

\section{Gauge equivalence of symplectic groupoids}

\label{sec:gaugegroup} In this section we will apply the results
of Section \ref{sec:predual} to dual pairs coming from symplectic
groupoids. 

A \textbf{symplectic groupoid} \cite{We87} is a symplectic manifold
\( (G,\Omega ) \) which is a Lie groupoid such that the graph 
\( \gamma _{m}=\{(x,y,m(x,y)),\; (x,y)\in G_{2}\} \)
of the multiplication is lagrangian in \( G\times G\times \overline{G} \)
(as usual, \( G_{2}\)
denotes the set of composable pairs). We denote the source (resp.
target) map of \( G \) by \( \alpha  \) (resp. \( \beta  \)), the
identity embedding by \( \epsilon :G_{0}\hookrightarrow G \), and
the inversion by \( i:G\longrightarrow G \). We recall that there
exists a unique Poisson structure \( \pi  \) on \( G_{0} \) making
\( \alpha  \) (resp. \( \beta  \)) into a Poisson (resp. anti-Poisson)
map.

Let \( (M,\pi ) \) be an integrable Poisson manifold, with symplectic
groupoid \( (G,\Omega ,\alpha ,\beta ) \). Since Lie algebroids corresponding
to gauge-equivalent Dirac structures are isomorphic, all the Lie algebroids
of Dirac structures in the \( \tau  \)-orbit of \( \pi  \) can be
integrated to a Lie groupoid isomorphic to \( (G,\alpha ,\beta ) \).
We now discuss the effect of a gauge transformation \( \tauB  \) on
the symplectic form \( \Omega  \). 

Let \( B \) be a closed \( 2 \)-form on \( M \), and consider the
\( 2 \)-form \( \Omega _{B} \doteq \Omega +\alpha ^{*}B-\beta ^{*}B \)
on \( G \). 

\begin{theorem} If \( \pi _{{B}}\doteq \tauB (\pi ) \) is Poisson, then
\( G_{{B}}\doteq (G,\Omega _{B},\alpha ,\beta ) \) is a symplectic groupoid
integrating \( (M,\pi _{{B}}) \). \end{theorem}

\begin{proof} By Theorem \ref{thm:ldualpair}, \( \Omega _{B} \)
is symplectic. We must check that the graph \[
\gamma _{m}=\{(x,y,m(x,y)),\; (x,y)\in G_{2}\}\]
 is lagrangian in \( G_{{B}}\times G_{{B}}\times \overline{G_{{B}}} \).
Let \( (x,y)\in G_{2} \), and consider a curve \( (x(t),y(t)) \)
in \( G_{2} \) with \( (x(0),y(0))=(x,y) \). Let \( (u,v)=(x'(0),y'(0)) \).
Then \( (u,v,T_{(x,y)}m(u,v))\in T_{p}\gamma _{m} \), \( p=(x,y,m(x,y)) \)
and any element in \( T_{p}\gamma _{m} \) is of this form. 

Differentiating the identities \( \alpha (m(x,y))=\alpha (x) \),
\( \beta (m(x,y))=\beta (y) \) and \( \beta (x)=\alpha (y) \), we
get 
\begin{gather}
T\alpha Tm (u,v) = T\alpha (u),\; \; T\beta Tm(u,v)= T\beta(v),
\label{eq:constraint1}\\
T\beta (u) = T\alpha (v).\label{eq:constraint2}
\end{gather} 
Therefore, if \( (u_{1},v_{1},Tm(u_{1},v_{1})),(u_{2},v_{2},Tm(u_{2},v_{2}))\in T_{p}\gamma _{m} \),
we have \begin{equation*}
\begin{split}
& (\Omega_{B}\times \Omega_{B} \times (-\Omega_{B}))
((u_{1},v_{1},Tm(u_{1},v_{1}),(u_{2},v_{2},Tm(u_{2},v_{2}))) = \\
& B(T\alpha (u_{1}),T\alpha (u_{2}))-B(T\beta (u_{1}),T\beta (u_{2})) + 
 B(T\alpha (v_{1}),T\alpha (v_{2}))-B(T\beta (v_{1}),T\beta (v_{2}))- \\ 
 & B(T\alpha Tm(u_{1},v_{1}),T\alpha Tm(u_{2},v_{2}) )+
 B(T\beta Tm(u_{1},v_{1}),T\beta Tm(u_{2},v_{2})) = 0,
\end{split}
\end{equation*}
Hence \( \gamma _{m} \) is lagrangian in \( G_{B}\times G_{B}\times \overline{{G}_{B}} \)
and \( G_{{B}} \) is a symplectic groupoid. 

By Theorem \ref{thm:gaugedual}, \( \alpha :(G,\Omega _{B})\longrightarrow (M,\piB ) \)
is a Poisson map. Since there is a unique Poisson structure on the
identity section of a symplectic groupoid with this property, the
Poisson structure induced by \( \Omega _{B} \) on \( M \) is \( \pi _{{B}} \).
\end{proof}

Thus, the effect of applying a gauge transformation \( \tauB  \)
to the Poisson structure of the identity section of a symplectic groupoid
is the following change of the symplectic form on the groupoid: \[
\Omega \stackrel{\tauB }{\longmapsto }\Omega _{B}=\Omega +\alpha ^{*}B-\beta ^{*}B.\]

Note that, in general,  \( \tauB (\pi ) \) is not Poisson, and the
form \( \Omega _{B} \) on \( G \) is degenerate. 
So 
\[\begin{diagram}
\node{}\node{(G,\Omega _{B})}\arrow{sw,t}{\alpha}\arrow{se,t}{\beta}\node{}\\
\node{(M,\piB)}\node{}\node{(M,-\piB)}
\end{diagram}\]
is generally just a pre-dual pair. In this context, one is naturally led to
consider groupoids equipped with pre-symplectic forms; 
the associated integration problem is whether a Dirac
structure with an integrable Lie algebroid can be integrated to such a 
``pre-symplectic groupoid'' (we believe an analog of the constructions
in \cite{CaFe,CrFe01} should clarify this question). Similarly, considering
gauge transformations on Poisson groupoids \cite{We88}, we are led
to the more general notion of ``Dirac groupoids''. 
The development of these ideas is the subject of work in progress.

\section{Morita equivalence of gauge-equivalent Poisson structures}

\label{sec:morita} In this section, we compare the notions of gauge and Morita equivalence
for integrable Poisson manifolds. 

Two Poisson manifolds \( (M_{1},\pi _{1}) \), \( (M_{2},\pi _{2}) \)
are called \textbf{Morita equivalent} \cite{Xu91} if there exists
a symplectic manifold \( (S,\Omega ) \) and Poisson maps \( J_{i}:S\longrightarrow M_{i} \),
\( i=1,2 \), so that 
\[\begin{diagram}
\node{}\node{(S,\Omega)}\arrow{sw,t}{J_1}\arrow{se,t}{J_2}\node{}\\
\node{(M_1,\pi_1)}\node{}\node{({M_2},-\pi_2)}
\end{diagram}\]
is a complete full dual pair with \( J_{i} \)-connected and \( J_{i} \)-simply-connected
fibers. In this case we call this diagram a  \textbf{Morita equivalence
bimodule}. 

Let \( (M,\pi ) \) be an integrable Poisson manifold with $\alpha$-connected
and $\alpha$-simply-connected symplectic groupoid \( (G,\Omega ,\alpha ,\beta ) \).
Let \( B \) be a closed \( 2 \)-form on \( M \) such that \( \piB =\tauB (\pi ) \)
is Poisson. 

\begin{theorem}\label{thm:morita} The Poisson manifolds \( (M,\pi ) \),
\( (M,\piB ) \) are Morita equivalent, with  Morita equivalence
bimodule \( (G,\widehat{\Omega },\alpha ,\beta ) \), 
where \( \widehat{\Omega }=\Omega -\beta ^{*}B \).
\end{theorem}

\begin{proof} By Theorem \ref{thm:gaugedual}, \( (M,\pi )\stackrel{\alpha }{\leftarrow }(G,\widehat{\Omega })\stackrel{\beta }{\rightarrow }(M,\tau _{-{B}}(-\pi )) \)
is a full dual pair with connected and simply connected fibers. Since
\( \tau _{{-B}}(-\pi ))=-\tauB (\pi ) \), \( \beta :(G,\widehat{\Omega })\longrightarrow (M,\tau _{{B}}(\pi )) \)
is anti-Poisson, and it only remains to show that this dual pair is
complete. 

\noindent Let \( X_{h} \) denote the hamiltonian vector field of
\( h \) with respect to \( \Omega  \). Let \( \widehat{X}_{h} \)
and \( X^{{B}}_{h} \) denote the hamiltonian vector fields with respect
to \( \widehat{\Omega } \) and \( \Omega _{B}=\Omega +\alpha ^{*}B-\beta ^{*}B \).

\noindent \textbf{Claim.} We have the following relations between
the hamiltonian vector fields:
$$
\widehat{X}_{{\alpha ^{*}f}}={X}_{{\alpha ^{*}f}} \;\; \mbox{ and }\;\;
\widehat{X}_{{\beta ^{*}f}}=X^{{B}}_{{\beta ^{*}f}}, \;\; \mbox{ for all }\; f \in C^\infty(M).
$$
\begin{proof} 
Since \( {X}_{{\alpha ^{*}f}}\in \ker \beta  \), \( \beta ^{*}B({X}_{{\alpha ^{*}f}})=0 \).
Hence \( \widehat{\Omega }({X}_{{\alpha ^{*}f}})=
{\Omega }({X}_{{\alpha ^{*}f}})=d(\alpha ^{*}f) \).
Therefore, \( \widehat{X}_{{\alpha ^{*}f}}={X}_{{\alpha ^{*}f}} \).
The other relation can be derived analogously. 
\end{proof}

Let now \( f\in C^{\infty }(M) \) be a complete function with respect
to \( \pi  \) (resp. \( \piB  \)). We must check that \( \widehat{X}_{\alpha ^{*}f} \)
(resp. \( \widehat{X}_{\beta ^{*}f} \)) is complete. 

By our claim, \( \widehat{X}_{\alpha ^{*}f}=X_{\alpha ^{*}f} \),
which is complete, since \( \alpha :(G,\Omega )\longrightarrow (M,\pi ) \)
is the target map of a symplectic groupoid \cite[Chp.~III]{CDW87}.
Analogously, since \( (G,\Omega _{B},\alpha ,\beta ) \) is a symplectic
groupoid for \( (M,\pi _{B}) \), it follows that \( \beta :(G,\Omega _{B})\longrightarrow (M,\piB ) \)
is complete, and hence \( \widehat{X}_{{\beta ^{*}f}}=X^{{B}}_{{\beta ^{*}f}} \)
is complete as well. \end{proof}

It is clear that Morita equivalent Poisson structures on a manifold
\( M \) need not be gauge equivalent, as their leaf decompositions
are not necessarily the same. More generally, we say that two Poisson structures
$\pi$ and $\pi'$ on $M$ are {\bf gauge equivalent up to Poisson diffeomorphism}
if there exists a diffeomorphism $f : M \longrightarrow M$
such that $f_*\pi$ and $\pi'$ are gauge equivalent.
It is clear that integrable Poisson structures which are
gauge equivalent up to Poisson diffeomorphism are still Morita equivalent.
However, as the next example
shows, one can have Morita equivalent Poisson structures on a manifold
\( M \) which are not gauge equivalent up to Poisson diffeomorphism. 

\begin{example} The example is based on non-cancellation properties
of product manifolds (see \cite{HilMisRoi71} and references therein).
Let \( F_{1},F_{2} \) and \( B \) be closed smooth manifolds so
that \( F_{1} \) and \( F_{2} \) have different homotopy types and
\( F_{1}\times B \) is diffeomorphic to \( F_{2}\times B \) (see
\cite{Charlap65} for original examples). As discussed in \cite{HilMisRoi71},
one can take \( B=S^{q} \) and \( F_{i} \) to be \( S^{q} \)-bundles
over \( S^{n} \), for \( q \) and \( n \) suitably chosen (and
large). Hence we can assume that \( B \) and \( F_{i} \) are simply
connected. 

Let \( E \) be the total space of these trivial fibrations, with
diffeomorphisms \( E\stackrel{\varphi _{i}}{\longrightarrow }{B\times F_{i}} \),
\( i=1,2 \). Let \( M=T^{*}E \). If \( \varphi _{i}^{\sharp } \)
denotes the natural cotangent lift of \( \varphi _{i} \), we obtain
the diffeomorphisms 
\[
\varphi _{i}^{\sharp }:M\longrightarrow T^{*}B\times T^{*}F_{i},\; \, i=1,2.
\]
 Since the transformations in the structure group of the fiber bundle
\[\begin{diagram}
\node{T^*F_i}\arrow{e,t}{}{}
\node{M}\arrow{s,r}{}\\
\node{}{}\node{T^*B}
\end{diagram}\]
are cotangent lifts
of diffeomorphisms of \( F_{i} \), they preserve the canonical
symplectic forms on \( T^*F_{i} \), \( i=1,2 \). Thus \( \varphi _{i}^{\sharp } \)
makes \( M \) into a bundle of symplectic manifolds in the sense
of \cite{GLSW83}, defining a Poisson structure \( \pi _{i} \) on
\( M \). As \( F_{1} \) and \( F_{2} \) have different homotopy
types, so do \( T^{*}F_{1} \) and \( T^{*}F_{2} \), and hence \( \pi _{1} \)
and \( \pi _{2} \) cannot be gauge equivalent up to Poisson diffeomorphism. 

Finally note that since \( F_{1},F_{2} \) and \( B \) are simply
connected, so are \( T^{*}F_{1},T^{*}F_{2} \) and \( T^{*}B \),
and \cite[Thm.~3]{GLSW83} implies that the symplectic structure along
the fibers of \( T^{*}F_{i}\rightarrow M\rightarrow T^{*}B \) admits
a closed extension. Hence the fundamental class \cite{DaDe87} of
\( (M,\pi _{i}) \) vanishes, and \cite[Thm.~4.3]{Xu91} implies that
\( (M,\pi _{i}) \) is Morita equivalent to the base \( T^{*}B \)
equipped with the zero Poisson structure. By transitivity, \( (M,\pi _{1}) \)
and \( (M,\pi _{2}) \) are Morita equivalent. 
\end{example}

\section{Gauge and Morita equivalence of
 topologically stable Poisson structures on surfaces}

\label{sec:applic}

Let \( \Sigma  \) be a compact connected oriented surface. Since
for dimensional reasons any bivector field on \( \Sigma  \) is Poisson,
Poisson structures on \( \Sigma  \) form a vector space. 

For \( n\geq 0 \), let \( \mathscr {G}_{n}(\Sigma ) \) be the set
of Poisson structures \( \pi  \) on \( \Sigma  \) such that 

\begin{itemize}
\item the zero set \(\{p\in \Sigma |\, \pi (p)=0\} \) consists
of \( n \) smooth disjoint curves \( \gamma _{1}(\pi ),\cdots ,\gamma _{n}(\pi ) \); 
\item \( \pi  \) vanishes linearly on each of the curves \( \gamma _{1}(\pi ),\cdots ,\gamma _{n}(\pi ) \). 
\end{itemize}
Let \( \mathscr {G}(\Sigma )\doteq \bigsqcup _{n\geq 0}\mathscr {G}_{n}(\Sigma ) \).
For \( n\geq 1 \) the symplectic leaves of \( \pi \in \mathscr {G}(\Sigma ) \)
are the points in the zero set \(\bigsqcup ^{n}_{i=1}\gamma _{i} \) and
the connected components of \( \Sigma \setminus \bigsqcup ^{n}_{i=1}\gamma _{i} \); for \( n=0 \)
the structure is symplectic. 
We call the Poisson structures in  \( \mathscr {G}(\Sigma )\) {\bf topologically stable},
since the topology of their zero sets is preserved under small perturbations.

Choosing a non-degenerate Poisson structure \( \pi _{0} \) on \( \Sigma  \),
we can identify the space of Poisson structures on \( \Sigma  \)
with \( C^{\infty }(\Sigma ) \): any \( \pi  \) is represented in
the form \( \pi =f\cdot \pi _{0} \) for \( f\in C^{\infty }(\Sigma ) \).
Under this identification, \( \mathscr {G}(\Sigma ) \) corresponds
to the space of smooth functions for which \( 0 \) is a regular value.
This implies that the set of topologically stable
Poisson structures \( \mathscr {G}(\Sigma ) \)
is \emph{generic}, i.e., \( \mathscr {G}(\Sigma ) \) is an open dense
subset of the space of all Poisson structures on \( \Sigma  \) endowed
with the Whitney \( C^{\infty } \) topology. 

As shown in \cite{Radko01}, the classification of Poisson structures
in \( \mathscr {G}(\Sigma ) \) up to Poisson isomorphisms depends
on a finite number of invariants. In order to recall what these invariants
are, we need a few definitions. 

For a Poisson manifold \( (M,\pi ) \), let \( \nu  \) be a volume
form on \( M \). The \textbf{modular vector field}~\( X^{\nu } \)
of \( \pi  \) with respect to \( \nu  \) \cite{We97} is defined
by the formula 
\[
X^{\nu }h\doteq\frac{L_{X_{h}}\nu }{\nu },\quad h\in C^{\infty }(M).
\]
This vector field measures the degree of invariance of \( \nu  \)
under the flows of hamiltonian vector fields; in particular, \( X^{\nu }=0 \)
if and only if \( L_{X_{h}}\nu =0 \) for all \( h\in C^{\infty }(M) \).
The modular vector field \( X^{\nu } \) is Poisson (i.e., its flow
preserves \( \pi  \)), and, if \( \nu ' \) is a different volume
form, one has
\[
X^{\nu '}=X^{\nu }+X_{-\log k},
\]
where \( k \) is the nowhere zero ratio \( k=|\frac{\nu '}{\nu }| \).
Thus the class of \( X^{\nu } \) in the first Poisson cohomology
\( H^{1}_{\pi }(M) \) (i.e., its equivalence class modulo hamiltonian
vector fields) is independent of \( \nu  \). This defines an element in
\( H^{1}_{\pi }(M) \)
called the \textbf{modular class} and denoted by \( \mu _{(M,\pi )} \). 

Suppose that the Poisson tensor \( \pi  \) on \( M \) vanishes on
a closed curve \( \gamma \subset M \), and is nonzero away from \( \gamma  \)
in a neighborhood of \( \gamma  \). Since the modular vector field
\( X^{\nu } \) preserves the Poisson structure, its flow must take
the zero set of \( \pi  \) to the zero set of \( \pi  \). Thus the
flow of \( X^{\nu } \) takes \( \gamma  \) to \( \gamma  \) and
so \( X^{\nu } \) must be tangent to \( \gamma  \). Moreover, for
another choice of volume form \( \nu ' \), we have 
\[
X^{\nu '}|_{\gamma }=X^{\nu }|_{\gamma }+(X_{-\log |\frac{\nu '}{\nu }|})|_{\gamma }=
X^{\nu }|_{\gamma },
\]
 since all hamiltonian vector fields are zero when restricted to the
zero curve \( \gamma  \). It follows that the restriction of the
modular vector field \( X^{\nu } \) to \( \gamma  \) is independent
of \( \nu  \), and hence induces an orientation on $\gamma$. 
As was observed in \cite{Roytenberg}, the period
of the flow of this vector field around \( \gamma  \) is an invariant
of the Poisson structure \( \pi  \). We denote this number by \( T_{\gamma }(M,\pi ) \)
(or, for short, \( T_{\gamma }(\pi ) \) when it is clear what \( M \)
is). 

For $\pi \in \mathscr {G}_{n}(\Sigma )$, let $Z(\pi )$ denote its zero set,
consisting on $n$ disjoint curves, taken with the induced orientations.
The main result of \cite{Radko01} states that Poisson structures
in \( \mathscr {G}_{n}(\Sigma ) \) are completely classified up to
(orientation-preserving) Poisson isomorphisms by the class of
\( Z(\pi )\) modulo orientation-preserving 
diffeomorphisms of $\Sigma$,
the \( n \) modular periods of \( \pi  \) around each connected
component of \( Z(\pi ) \) and the regularized Liouville volume $V(\pi)$(which
is a certain regularized sum of symplectic volumes of two-dimensional
leaves, taken with appropriate signs). 

We shall now consider the questions
of gauge and Morita equivalence of Poisson structures in \( \mathscr {G}(\Sigma ) \).

\subsection{Gauge equivalence in \protect\( \mathscr {G}(\Sigma )\protect \)}

\label{subsec:applicgauge}

The obvious necessary condition for two Poisson structures \( \pi  \),
\( \pi '\in \mathscr {G}(\Sigma ) \) to be gauge equivalent is \( Z(\pi )=Z(\pi ') \),
i.e. the zero sets of both structures, with the induced orientations, should be the same. 

\begin{proposition}\label{gauge_pres_periods}
Let  \( \pi ,\, \pi '\in \mathscr {G}_{n}(\Sigma ) \) be Poisson structures with
\( Z(\pi )=Z(\pi ')=\bigsqcup _{i=1}^{n}\gamma _{i} \).
If they are gauge equivalent, then their modular periods are the
same around all  the zero curves, i.e. \( T_{\gamma _{i}}(\pi )=T_{\gamma _{i}}(\pi ') \),
for \( i=1,\dots ,n \). 
\end{proposition}

\begin{proof}
Let \( \pi =f\cdot \pi _{0} \), \( \pi '=f'\cdot \pi _{0} \),
where \( f,\, f'\in C^{\infty }(\Sigma ) \) are functions vanishing
linearly on \( \gamma _{1},\cdots ,\gamma _{n} \) and non-zero elsewhere,
and \( \pi _{0} \) is a fixed nondegenerate Poisson structure on
\( \Sigma  \). 

For each \( i=1,\dots ,n \), let 
\( U_{i}=\{(z_{i},\theta _{i})|\, |z_{i}|<R_{i},\, \theta _{i}\in [0,2\pi ]\} \)
be a small annular neighborhood of the zero curve \( \gamma _{i}\in Z(\pi ) \)
such that \( Z(\pi )\cap U_{i}=\gamma _{i} \) and 
\( \pi |_{U_{i}}=f(z_{i},\theta_i)\partial _{z_{i}}\wedge \partial _{\theta _{i}} \),
\( \pi '|_{U_{i}}=f'(z_{i},\theta_i)\partial _{z_{i}}\wedge \partial _{\theta _{i}} \).
A simple computation shows that the modular vector fields along $\gamma_i$ are
$$
X_{\pi}= \frac{\partial f}{\partial z_i}(0,\theta_i)\partial_{\theta_i},\;\;\;
X_{\pi'}= \frac{\partial f'}{\partial z_i}(0,\theta_i)\partial_{\theta_i}.
$$
Suppose that \( \pi '=\tauB(\pi ) \) for a closed \( 2 \)-form
\( B\in \Omega ^{2}(\Sigma ) \). Writing $B = b(z_i,\theta_i)dz_i \wedge d\theta_i$, it follows that
$$
f'(z_i,\theta_i)= \frac{f(z_i,\theta_i)}{1 + f(z_i,\theta_i)b(z_i,\theta_i)}, 
$$
which implies that
$$
\frac{\partial f}{\partial z_i}(0,\theta_i)=\frac{\partial f'}{\partial z_i}(0,\theta_i).
$$
So $\pi$ and $\pi'$ have equal modular vector fields along $\gamma_i$, and hence equal modular
periods.
\end{proof}

We remark that if
$\pi \in \mathscr {G}_{n}(\Sigma)$ and $B$ is a closed $2$-form on $\Sigma$
such that $\tauB(\pi)$ is Poisson, then the regularized Liouville volumes
of $\pi$ and $\tauB(\pi)$ are related by
\( V(\tauB(\pi))=V(\pi )+\textrm{Vol}(B) \), where
\( \textrm{Vol}(B)=\int _{\Sigma }B \) is the Liouville volume of
\( B \).

We now discuss the converse of Proposition \ref{gauge_pres_periods}.

\begin{theorem} \label{Thm:gauge_equiv_necklaces} 
Let \( \pi ,\, \pi '\in \mathscr {G}_{n}(\Sigma ) \) be two Poisson structures
with \( Z(\pi )=Z(\pi ')=\bigsqcup _{i=1}^{n}\gamma _{i} \).
If \( T_{\gamma _{i}}(\pi )=T_{\gamma _{i}}(\pi ') \),
for \( i=1,\dots ,n \), then $\pi$ and $\pi'$ are gauge equivalent up to Poisson
diffeomorphism.
\end{theorem}

\begin{proof} 
Let \( \pi =f\cdot \pi _{0} \), \( \pi '=f'\cdot \pi _{0} \), 
with $f$ and $f'$ as in the proof of Prop.~\ref{gauge_pres_periods}.
By replacing, if necessary,
 $\pi'$ by a Poisson diffeomorphic structure (with same zero set),
we can assume that, for each $i=1,\dots,n$, there is a small annular
neighborhood 
\( U_{i}=\{(z_{i},\theta _{i})|\, |z_{i}|<R_{i},\, \theta _{i}\in [0,2\pi ]\} \)
of  \( \gamma _{i}\in Z(\pi ) \) 
such that \( Z(\pi )\cap U_{i}=\gamma _{i} \)
and \( \pi |_{U_{i}}= f(z_{i})\partial _{z_{i}}\wedge \partial _{\theta _{i}} \),
\( \pi '|_{U_{i}}=f'(z_{i})\partial _{z_{i}}\wedge \partial _{\theta _{i}} \)
with $ f|_{U_{i}}=c_{i}z_{i}+O(z^{2})$, and $ f'|_{U_{i}}=c'_{i}z_{i}+O(z^{2})$.
A simple computation (see (\ref{eq:modvect})) shows that 
$c_i = \frac{2\pi}{T_{\gamma _{i}}(\pi)}$ and 
$c'_i = \frac{2\pi}{T_{\gamma _{i}}(\pi')}$, so $c_i=c_i'$.

Let $\omega = f^{-1}\omega_0$ (resp. $\omega' = f'^{-1}\omega_0$) be the symplectic
form corresponding to $\pi$ (resp. $\pi'$) on $\Sigma \setminus Z(\pi)$, and
define the $2$-form $B$ on $\Sigma \setminus Z(\pi)$ by
$$
B=\omega' - \omega = \left( \frac{1}{f'}-\frac{1}{f}\right)\omega_0.
$$
It is simple to check that, since $c_i = c_i'$, the function
$\frac{1}{f'} - \frac{1}{f} = \frac{f-f'}{f f'}$ extends to a smooth function on $\Sigma$.
Hence the $2$-form $B$ can be extended to a (closed) $2$-form
on $\Sigma$, also denoted by $B$, with the property that
$\pi' = \tauB(\pi)$.
\end{proof}

For \( \pi \in \mathscr {G}_{n}(\Sigma ) \), 
the anchor \( \tilde{\pi }:\, T^{*}\Sigma \rightarrow T\Sigma  \)
of the corresponding
Lie algebroid is injective on the open dense set \( \Sigma \setminus Z(\pi ) \).
According to \cite[Thm.~1]{Debord00}, a Lie algebroid whose anchor
is injective on an open dense set is integrable, so \( (\Sigma ,\pi ) \)
is an integrable Poisson manifold. 

The following result follows from Theorem \ref{thm:morita}.

\begin{theorem} Two Poisson structures \( \pi ,\, \pi '\in \mathscr {G}_{n}(\Sigma ) \)
with the same zero sets \( Z(\pi )=Z(\pi ')=\bigsqcup _{i=1}^{n}\gamma _{i} \)
and equal modular periods, \( T_{\gamma _{i}}(\pi )=T_{\gamma _{i}}(\pi ') \),
for \( i=1,\dots ,n \), are Morita equivalent. 
\end{theorem}

We now turn our attention to the study of Morita equivalence in \( \mathscr {G}(S^{2}) \).
First, we need to collect a few general results on 
invariants of Morita equivalence.

\subsection{Invariants of Morita equivalence: topology of the leaf space and modular
periods}

\label{subsec:applicmorita}

Let \( (M,\pi ) \) be a Poisson manifold.
Let $L(M)$ be the leaf space of the symplectic foliation of $\pi$,
endowed with its quotient topology:
for a topological space \( X \), a function \( f:L(M)\rightarrow X \)
is continuous if and only if \( f\circ \pr :M\rightarrow X \) is continuous,
where  \( \pr :M\rightarrow L(M) \)
is the quotient map. 

Let \( (M_{1},\pi _{1}) \) and \( (M_{2},\pi _{2}) \) be Poisson
manifolds, and let \( (M_{1},\pi _{1})\stackrel{J_{1}}{\leftarrow }(S,\Omega )\stackrel{J_{2}}{\rightarrow }(M_{2},\pi _{2}) \)
be a Morita equivalence bimodule. It is well-known (see e.g. \cite{Blaom,SilWein99})
that \( S \) induces a bijection of sets \( \phi _{S}:L(M_{1})\rightarrow L(M_{2}) \)
given by
\[
\phi _{S}(\mathscr {L})=J_{2}(J_{1}^{-1}(\mathscr {L})),\quad 
\textrm{for }\mathscr {L}\in L(M_{1}).\]
The following observation is based on ideas from \cite{Crainic}. 

\begin{proposition} \label{Prop:Morita_inv_topology}The map \( \phi _{S}:L(M_{1})\rightarrow L(M_{2}) \)
is a homeomorphism of topological spaces. \end{proposition}

\begin{proof} Let \( F_{i} \) be the subset of \( TM_{i} \) consisting
of vectors tangent to the symplectic leaves. Let \( TJ_{i}\subset TS \),
\( i=1,2 \) be the subbundles tangent to the \( J_{i} \)-fibers.
Then\[
J_{1}^{*}F_{1}=J_{2}^{*}F_{2}=TJ_{1}+TJ_{2},\]
 where \( J_{i}^{*}F_{i}=\{v\in TS|\, TJ_{i}v\in F_{i}\} \) denotes
the pull-back of \( F_{i} \). Let \( F=J_{1}^{*}F_{1}=J_{2}^{*}F_{2} \). 

Since the fibers of \( J_{i} \)  are connected, \( i=1,2 \), the
natural maps 
\[
\psi _{i}:S/F\rightarrow M_{i}/F_{i}\]
 of leaf spaces are bijections. Moreover, it is not hard to see that
\( \phi _{S}  =\psi _{2}\circ \psi _{1}^{-1} \).
So, if we  endow \( S/F \) with its quotient topology,  it is sufficient
to prove that \( \psi _{i} \), \( i=1,2 \), are homeomorphisms. 

By the definition of the quotient topology, the map \( \psi _{i}:S/F\rightarrow M_{i}/F_{i} \)
is continuous if and only if the map \( \psi _{i}\circ \pr :S\rightarrow M_{i}/F_{i} \)
is continuous (here \( \pr :S\rightarrow S/F \) is the quotient map).
But \( \psi _{i}=\pr _{i}\circ J_{i} \), where \( \pr _{i}:M_{i}\rightarrow M_{i}/F_{i} \)
is the quotient map. Hence \( \psi _{i} \) is continuous. 

Similarly, \( \psi _{i}^{-1}:M_{i}/F_{i}\rightarrow S/F \) is continuous
if and only if \( \psi _{i}^{-1}\circ \pr _{i}:M_{i}\rightarrow S/F \) is continuous.
Since \( J_{i} \) is a submersion, this is true if and only if \( \psi _{i}^{-1}\circ \pr _{i}\circ J_{i}:S\rightarrow S/F \)
is continuous. But \( \psi _{i}^{-1}\circ \pr _{i}\circ J_{i}=\pr  \).
Therefore, \( \psi _{i} \) is a homeomorphism for \( i=1,2 \), which
implies that \( \phi _{S} \) is a homeomorphism. 
\end{proof}

The modular vector field and the modular class are well-behaved under
Morita equivalence: the bijection of leaf spaces \( \phi _{S} \)
induced by a Morita equivalence bimodule 
produces an isomorphism of Poisson cohomologies \cite[Thm.~3.1]{GinzLu92}
\begin{equation}
\label{isomorphism_of_H1}
\phi _{S}^{*}:H^{1}_{\pi }(M_{1})\rightarrow H^{1}_{\pi }(M_{2}),
\end{equation}
 which preserves the modular class \cite{Crainic,Ginz2000}, \[
\phi ^{*}_{S}(\mu _{(M_{1},\pi _{1})})=\mu _{(M_{2},\pi _{2})}.\]
 We will need the following remark from the construction of the isomorphism
(\ref{isomorphism_of_H1}) in \cite{GinzLu92}. 

\begin{remark} \label{rmq:existsModularOnS}
Given volume forms \( \nu _{1} \)
and \( \nu _{2} \) on \( M_{1} \) and \( M_{2} \), respectively,
there exists a vector field \( X \) on \( S \) with the property
that \( (J_{i})_{*}X=X^{\nu _{i}} \), \( i=1,2 \). The vector field
\( X \) is actually hamiltonian, and its Hamiltonian \( H \) is
determined by the equation\[
DJ_{1}^{*}\nu _{1}=\pm e^{H}J_{2}^{*}\nu _{2},\]
 where \( D:\Omega ^{k}(M)\rightarrow \Omega ^{2m-k}(M) \) is the
symplectic \( * \)-operator (cf. \cite{Bryl88}). \end{remark}

We now observe that the modular periods also behave well under Morita
equivalence. 

\begin{theorem} Let \( (M_{i},\pi _{i}) \), \( i=1,2 \) be Poisson
manifolds, and let \( (M_{1},\pi _{1})\stackrel{J_{1}}{\leftarrow }(S,\Omega )\stackrel{J_{2}}{\rightarrow }(M_{2},\pi _{2}) \)
be a Morita equivalence bimodule. Assume that \( Z_{i}\subset M_{i} \)
are such that \( \pi _{i}|_{Z_{i}}=0 \) and the isomorphism of leaf
spaces satisfies\[
\phi _{S}(Z_{1})=Z_{2}.\]
 Let \( \Phi ^{i}_{t} \) be the flow of the modular vector field
\( X^{\nu _{i}} \) for some volume form \( \nu _{i} \) on \( M_{i} \),
\( i=1,2 \). Assume that \( \Phi _{t}^{i} \) takes \( Z_{i} \)
to \( Z_{i} \) for all \( t \). Then\[
\phi _{S}\circ \Phi _{t}^{1}=\Phi _{t}^{2}\circ \phi _{S}~\mbox{ in } Z_1,\quad \forall t\in \mathbb {R}.\]
\end{theorem}

\begin{proof} Let \( X \) be a vector field on \( S \) such that
\( (J_{i})_{*}X=X_{i}^{\nu _{i}} \), \( i=1,2 \) (see Remark \ref{rmq:existsModularOnS}).
By the definition of \( \phi _{S} \), for each \( p_1\in Z_{1} \)
(which forms by itself a symplectic leaf), \[
\phi _{S}(\{p_1\})=J_{2}(J_{1}^{-1}(\{p_1\})),\]
 which by our assumption on \( \phi _{S} \) is a single point \( p_2\in Z_{2} \).
It follows that \( J_{1}^{-1}(\{p_1\})\subset J_{2}^{-1}(\{p_2\}) \).
Reversing the roles of \( p_1 \) and \( p_2 \), we get 
\( J_{1}^{-1}(\{p_1\})=J_{2}^{-1}(\{p_2\}). \)
In particular, it follows that \[
J_{1}^{-1}(Z_{1})=J_{2}^{-1}(Z_{2}).\]
 Thus for any fixed \( r\in J_{1}^{-1}(\{p_1\})=J_{2}^{-1}(\{p_2\}) \)
we obtain\[
\Phi _{t}^{2}(\phi _{S}(\{p_1\}))=J_{2}(\Phi _{t}(\{r\}))=
J_{2}(J_{1}^{-1}(\Phi ^{1}_{t}(\{p_1\}))=\phi _{S}(\Phi _{t}^{1}(\{p_1\})),
\]
where $\Phi _{t}$ is the flow of $X$.
 Therefore, \( \phi _{S}\circ \Phi _{t}^{1}=\Phi _{t}^{2}\circ \phi _{S} \) in $Z_1$.
\end{proof}

\begin{corollary} \label{Morita_inv_of_modular_periods}Assume that
\( \gamma _{i}\subset M_{i} \) are simple closed curves, \( U_{i}\supset \gamma _{i} \)
are open sets, and that \( \pi _{i}|_{\gamma _{i}}=0 \) and \( \pi |_{U_{i}\setminus \gamma _{i}}\neq 0 \),
\( i=1,2 \). Assume further that \( (M,\pi _{1}) \) and \( (M_{2},\pi _{2}) \)
are Morita-equivalent via a bimodule \( S \). Assume finally that
\( \phi _{S}(\gamma _{1})=\gamma _{2} \). If the modular vector fields are nonzero in $\gamma_i$,$i=1,2$,
then \[T_{\gamma _{1}}(M_{1},\pi _{1})=T_{\gamma _{2}}(M_{2},\pi _{2}).\]
\end{corollary}

\begin{proof} Applying the previous theorem with \( Z_{i}=\gamma _{i} \),
\( i=1,2 \), we obtain that the restriction to $\gamma_i$ of the flows \( \Phi ^{i}_{t} \) of modular
vector fields \( X_{i}^{\nu _{i}} \) are intertwined by \( \phi _{S} \).
Thus, for any \( p_1\in \gamma _{1} \), 
\begin{eqnarray*}
T_{\gamma _{1}}(M_{1},\pi _{1}) & = & \inf \{t>0:\Phi ^{1}_{t}(p_1)=p_1\}\\
 & = & \inf \{t>0:\phi _{S}(\Phi ^{1}_{t}(p_1))=\phi _{S}(p_1)\}\\
 & = & \inf \{t>0:\Phi _{t}^{2}(\phi _{S}(p_1))=\phi _{S}(p_1)\}\\
 & = & \inf \{t>0:\Phi _{t}^{2}(p_2)=p_2\},\qquad p_2=\phi _{S}(p_1)\\
 & = & T_{\gamma _{2}}(M_{2},\pi _{2}).
\end{eqnarray*}
\end{proof}

\subsection{Morita equivalence of topologically stable Poisson 
structures on \protect\protect\( S^{2}\protect \protect \)}
\label{subsec:applicS2} 
We will now show that two 
Poisson structures in \( \mathscr {G}_{n}(S^{2}) \) 
are Morita equivalent if and only if they are gauge equivalent
up to diffeomorphism.
As a result, we find a complete Morita-equivalence invariant for
topologically stable Poisson structures on $S^2$.

Let \( \Sigma=S^{2} \), the \( 2 \)-sphere. Let \( \pi \in \mathscr {G}(\Sigma) \)
and write, as before, \( \pi =f\cdot \pi _{0} \), where \( \pi _{0} \)
is a non-degenerate Poisson structure on \( \Sigma \) and \( f\in C^{\infty }(\Sigma) \).
We fix on $\Sigma$ the orientation induced by (the symplectic form corresponding to)
$\pi_0$.
The class of \( Z(\pi ) \) modulo (orientation-preserving)
diffeomorphisms of \( \Sigma \) 
can be described combinatorially by a signed tree
in the following way \cite[Sec.~2.8]{Radko01}. The vertices
of the tree correspond to the connected components of \( \Sigma\setminus Z(\pi ) \)
(i.e., the regions bounded by the curves comprising \( Z(\pi ) \)).
The vertex is assigned a positive sign if \( f \) is positive on
the corresponding region; otherwise the vertex is assigned a negative
sign. Two vertices are connected
by an edge if and only if the corresponding regions share a boundary.
We denote this signed tree by \( Tree(\pi ) \). 

\begin{lemma} \label{Lemma:neccCondMoritaEquiv}If two  Poisson
structures \( \pi ,\pi '\in \mathscr {G}_{n}(S^2) \) are Morita equivalent,
then there exists an isomorphism of trees \( \phi :Tree(\pi )\rightarrow Tree(\pi ') \),
not necessarily preserving signs,
so that \( T_{\gamma }(\pi )=T_{\phi (\gamma )}(\pi ') \) for every
edge \( \gamma  \) of \( Tree(\pi ) \). 
\end{lemma}

\begin{proof} Assume that \( (S^2,\pi ) \) and \( (S^2,\pi ') \) are
Morita equivalent. Let \( \phi  \) be the induced homeomorphism of
the leaf spaces of \( (S^2,\pi ) \) and \( (S^2,\pi ') \) (see Prop.~\ref{Prop:Morita_inv_topology}). 

As a set, the leaf space \( L \) of \( (S^2,\pi ) \) can be identified
with the union \( Z(\pi )\sqcup \{\ell _{1},\cdots ,\ell _{n}\} \),
where \( \ell _{1},\dots ,\ell _{n} \) are the points corresponding
to the \( 2 \)-dimensional leaves \( \mathcal{L}_{1},\dots ,\mathcal{L}_{n} \).
The quotient topology of \( L \) is easily described: the only open
subsets of \( L \) have the form \( U\cup \{\ell _{i_{1}}\}\cup \cdots \cup \{\ell _{i_{k}}\} \),
where \( i_{1},\dots ,i_{k}\in \{1,\cdots, n\} \), \( k\geq 0 \) and \( U\subset Z(\pi ) \)
is an open subset with the property that if \( U \) intersects non-trivially
a curve \( \gamma \subset Z(\pi ) \), then for both leaves bounding
\( \gamma  \) the corresponding points of the leaf space occur among
\( \{\ell _{i_{1}},\ldots ,\ell _{i_{k}}\} \). 

Given \( L \) with its topology, consider the collection \( \mathcal{Y} \)
of all subsets \( Y\subset L \) with the property that \( L\setminus {Y} \)
is Hausdorff. Order \( \mathcal{Y} \) by inclusion. We claim that
\( X=\{\ell _{1},\ldots ,\ell _{n}\} \) is a minimal element of \( \mathcal{Y} \)
of smallest cardinality. First, note that \( X\in \mathcal{Y} \),
since the relative topology on \( L\setminus X=Z(\pi )\subset L \)
is Hausdorff. Next, assume that \( Y\in \mathcal{Y} \), and \( Y\not \supset X \).
Then \( \ell _{i}\in L\setminus Y \) for some \( i\in \{1,\cdots,n\} \).
Now all of the points of the boundary of \( \mathcal{L}_{i} \) in
\( S^2 \) lie in \( Z(\pi ) \) and cannot be separated from \( \ell _{i} \)
by open sets; thus all of these points must necessarily be in \( Y \).
Thus \( Y \) must have infinite cardinality. 

It follows that \( \phi  \) must map \( X \) to a subset of \( L' \)
with the same minimality property; and hence \( \phi  \) must take
the complement of \( X \), \( Z(\pi ) \), to \( Z(\pi ') \). Thus
\( \phi  \) induces a map between the set of vertices of \( Tree(\pi ) \)
and \( Tree(\pi ') \). 

Now, two vertices \( \ell _{i},\ell _{j}\in Tree(\pi ) \) are connected
by an edge if and only if the corresponding regions share a boundary in \( S^2 \).
A point \( x\in Z(\pi )\subset L \) cannot be separated from \( \ell _{k} \)
by an open set if and only if \( x \) belongs to the boundary of
\( \mathcal{L}_{k} \) in \( S^2 \). It follows that \( \ell _{i},\ell _{j} \)
are connected by an edge if and only if there exists a point \( x\in L \), such
that \( x\neq \ell _{i} \), \( x\neq \ell _{j} \), but which cannot
be separated from either of them by an open set. Since \( \phi  \)
is a homeomorphism, it must preserve this property, and thus \( \phi  \)
induces a map of trees from \( Tree(\pi ) \) to \( Tree(\pi ') \).
The statement about modular periods now follows from Corollary \ref{Morita_inv_of_modular_periods}.
\end{proof}

\begin{theorem} Two Poisson structures \( \pi ,\pi '\in \mathscr {G}_{n}(S^{2}) \)
are Morita equivalent if and only if there exists an isomorphism of
trees \( \phi :Tree(\pi )\rightarrow Tree(\pi ') \), not necessarily
preserving signs, so that \( T_{\gamma }(\pi )=T_{\phi (\gamma )}(\pi ') \)
for every edge \( \gamma  \) of \( Tree(\pi ) \). \end{theorem}

\begin{proof} The necessity follows from Lemma \ref{Lemma:neccCondMoritaEquiv}. 

Assume now that there exists an isomorphism \( \phi :Tree(\pi )\rightarrow Tree(\pi ') \)
satisfying the conditions. Let \( \psi :S^{2}\rightarrow S^{2} \)
be an orientation-reversing diffeomorphism. By replacing \( \pi ' \)
with \( \psi _{*}\pi ' \) (which is obviously Poisson isomorphic,
and hence Morita equivalent to, \( \pi ' \)) if necessary, we may
assume that \( \phi :Tree(\pi )\rightarrow Tree(\pi ') \) is an isomorphism
of \emph{signed} trees. 

Choose a function \( g\in C^{\infty }(S^{2}) \) supported on the
interior of one of the two-dimensional leaves. Let \( \pi ''=\pi '+g\pi  \).
Since \( \pi ''=\pi ' \) in a neighborhood of each of the zero curves
\( \gamma \subset Z(\pi ) \), the modular periods of \( \pi ' \)
and \( \pi '' \) are equal. Therefore, by Theorem \ref{Thm:gauge_equiv_necklaces},
\( \pi ' \) and \( \pi '' \) are gauge equivalent up to Poisson diffeomorphism, and hence Morita
equivalent (Theorem~\ref{thm:morita}) for any such choice of \( g \).
Also, the isomorphism \( \phi  \) induces an isomorphism of trees
\( \phi ':Tree(\pi )\rightarrow Tree(\pi '') \). 

With a suitable choice of \( g \), the regularized Liouville volume
of \( \pi '' \) can be made equal to that of \( \pi  \) (see \cite{Radko01}
for details). Thus by \cite[Thm.~3]{Radko01}, \( \pi  \)
and \( \pi '' \) are Poisson isomorphic. We conclude that \( \pi  \)
and \( \pi '' \) are Morita equivalent, and so are \( \pi  \) and
\( \pi ' \), by transitivity of Morita equivalence. \end{proof}

\begin{footnotesize}


\end{footnotesize}
\end{document}